\documentclass[oneside]{amsart}
\usepackage{amsmath, amssymb,color,amsthm,graphicx,bm,cite}
\usepackage{color,bm}
\newtheorem{theorem}{Theorem}[section]
\newtheorem{proposition}[theorem]{Proposition}
\newtheorem{lemma}[theorem]{Lemma}
\newtheorem{corollary}[theorem]{Corollary}

\newtheorem{question}{Question}

\theoremstyle{definition}
\newtheorem{definition}{Definition}
\newtheorem{main}{Theorem}

\def\Z{\mathbb{Z} }
\def\R{\mathbb{R} }

\def\nbd{neighborhood }
\def\nbds{neighborhoods }
\def\SX{\mathop{\mathrm{Sing}}(X)}
\def\PX{\mathop{\mathrm{Per}}(X)}
\def\CX{\mathop{\mathrm{Cl}}(X)}

\author{Tomoo Yokoyama}
\date{\today}
\address{Department of Mathematics, Faculty of Science, Saitama University, Shimo-Okubo 255, Sakura-ku, Saitama-shi, 338-8570 Japan\\}
\email{tyokoyama@rimath.saitama-u.ac.jp}
\thanks{The author was partially supported by JSPS Grant Numbers 20K03583}

\subjclass[2010]{}

\title[Combinatorial structures of the space of Hamiltonian vector fields]{Combinatorial structures of the space of Hamiltonian vector fields on compact surfaces}
\makeatletter
\@namedef{subjclassname@2020}{%
  \textup{2020} Mathematics Subject Classification}
\makeatother


\begin{document}
\maketitle

\begin{abstract}
In the time evolution of fluids, the topologies of fluids can be changed by the creations and annihilations of singular points and by switching combinatorial structures of separatrices. In this paper, to describe the possible generic time evolution of Hamiltonian vector fields on surfaces with or without constraints, we study the structure of the ``moduli space'' of such vector fields under the non-existence of creations and annihilations of singular points. In fact, we describe the relations of bifurcations between Hamiltonian vector fields to construct foundations of descriptions of the time evaluations of fluid phenomena. Moreover, we show that the space of topologically equivalence classes of such vector fields has non-contractible connected components and is a disjoint union of finite abstract cell complexes such that the codimension of a cell corresponds to the instability of a Hamiltonian vector field by using combinatorics and simple homotopy theory. In particular, there is a connected component of the space that is weakly homotopic to a three-dimensional sphere.
\end{abstract}

\section{Introduction}\label{intro}

\begin{figure}[t]
\begin{center}
\includegraphics[scale=0.45]{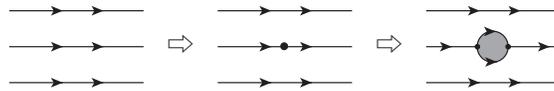}
\end{center}
\caption{Creation of a physical boundary.}
\label{fig:creations_bdry}
\end{figure}

In the time evolution of incompressible fluids on punctured spheres, the topologies of such fluids can be changed by the creations and annihilations of singular points and physical boundaries. 
For instance, one can observe the creation of a physical boundary, which is a boundary of a stone on the surface of a river, when the water level of the river goes down as in Figure~\ref{fig:creations_bdry}. 
Notice that creations and annihilations of physical boundaries change the topologies of surfaces. 
In other words, no creations and annihilations of physical boundaries appear when we consider flows on a fixed surface. 
On the other hand, the topologies of such fluids also can be changed by switching combinatorial structures of separatrices. 
Such combinatorial structures are studied from fluid mechanics, integrable systems, and dynamical systems. 
For instance, the topological classification of Hamiltonian flows with finitely many singular points has been investigated on a plane $\R^2$ \cite{aref1998stagnation}, on a sphere $\mathbb{S}^2$ \cite{kidambi2000streamline}, and on a torus \cite{moffatt2001topology} from a fluid mechanics point of view. 
These classifications are generalized to such flows on closed surfaces from an integrable system point of view \cite{bolsinov2004integrable} and compact surfaces from a dynamical system point of view \cite{sakajo2014unique,sakajo2015transitions,sakajo2018tree,sakajo2020discrete,yokoyama2021complete,yokoyama2013word}. 
%
In the codimension zero and one cases, such graph structures are characterized \cite{ma2005geometric,sakajo2015transitions,sakajo2018tree,yokoyama2013word,yokoyama2021cot}. 
In particular, vertices correspond to structurally stable Hamiltonian vector fields and edges correspond to ``generic'' intermediate vector fields.  
However, ``non-generic'' intermediate vector fields (e.g. fluids with symmetric vortex pairs) naturally appear under physical or symmetric restrictions. 
For instance, degenerate multi-saddle connections are the reasons for ``non-genericity''. 
Such degenerate multi-saddle connections for Hamiltonian vector fields on closed surfaces are studied from the integrable system's point of view \cite{bolsinov2004integrable}. 
More generally, we would like to ask whether the structure of the ``moduli space'' of Hamiltonian vector fields under the non-existence of creations and annihilations of singular points is ``nontrivial''. 
More precisely, we ask the following question. 
\begin{question}
Does the space of topological equivalence classes of Hamiltonian vector fields on a compact surface have non-contractible connected components, under the non-existence of creations and annihilations of singular points? 
\end{question}

In this paper, we demonstrate that there is such a non-contractible connected component as follows. 

\begin{main}\label{main:01-}
For any $r \in \Z_{\geq 0} \sqcup \{ \infty \}$, the space of $C^r$ Hamiltonian vector fields with finitely many singular points and with exactly $3$ topological centers but without fake multi-saddles \rm{(i.e.} $\partial$-$0$-saddles and $0$-saddles \rm{)} on a closed disk has a connected component whose weak homotopy type is a three-dimensional sphere. 
\end{main}

This means that the structure of ``transitions of Hamiltonian vector fields on compact surfaces'' is not contractible. 
To calculate the non-contractible property, we study the hierarchical structure of the space of such vector fields under the non-existence of creations and annihilations of singular points (see Theorem~\ref{main:01} and Theorem~\ref{main:02} for details). 
In fact, the space of topologically equivalence classes of such vector fields is a disjoint union of finite abstract cell complexes such that the codimension of a cell corresponds to the instability of a Hamiltonian vector field as follows. 


\begin{main}\label{main:02-}
The following statements hold for any $l \in \Z_{\geq 0}$ and any $r \in \Z_{\geq 0} \sqcup \{ \infty \}$:
\\
{\rm(1)} The space $X$ of topologically equivalence classes of $C^r$ Hamiltonian vector fields with at most $l$ singular points but without fake multi-saddles on a compact surface is a finite abstract cell complex, a finite poset, and a finite $T_0$-space. 
\\
{\rm(2)} For any $h \in \Z_{\geq 0}$, the closure of the set of elements of height $h$ in $X$ is the set  of elements of height at most $h$ in $X$ \rm{(i.e.} $\overline{X_{=h}} = \bigsqcup_{j=0}^h X_{=j}$, where $X_{=s}$ is the set of elements of height $s$ in the poset $X$ \rm{)}. 
\end{main}

Using the hierarchical structure, we can describe the possible generic time evolution of Hamiltonian vector fields on a compact surface with or without restriction conditions.

%

The present paper consists of six sections.
In the next section, as preliminaries, we introduce fundamental concepts.
In \S 3, we demonstrate the properties of perturbations of Hamiltonian vector fields on surfaces. 
In \S 4, we construct the filtration and abstract cell complex structures of the space of Hamiltonian vector fields on compact surfaces. 
In \S 5, the non-contractible property is demonstrated. 
In the final section, the application of this research method is described.
In particular, we explain that a similar approach can be applied with respect to gradient flow. 

\section{Preliminaries}\label{sec:prel}

We recall notions of order, topology, and dynamical systems to describe the main results precisely. 

%

\subsection{Notion of orders}

A binary relation $\leq$ on a set $X$ is a {\bf pre-order} if it is reflexive (i.e. $a \leq a$ for any $a \in X$) and transitive (i.e. $a \leq c$ for any $a, b, c \in X$ with $a \leq b$ and $b \leq c$).
For a pre-order $\leq$, the inequality $a<b$ means both $a \leq b$ and $a \neq b$.
A pre-order $\leq$ on $X$ is a partial order if it is antisymmetric (i.e. $a = b$ for any $a,b \in X$ with $a \leq b$ and $b \leq a$).
A poset is a set with a partial order.
A pre-order order $\leq$ is a total order (or linear order) if either $a < b$ or $b < a$ for any  points $a \neq b$.
A {\bf chain} is a totally ordered subset of a pre-ordered set with respect to the induced order.
Let $(X, \leq)$ be a pre-ordered set.
For a point $x \in X$, define the {\bf upset} $\mathop{\uparrow} x  := \{ y \in X \mid x \leq y \}$ of $x$ and the {\bf downset} $\mathop{\downarrow} x  := \{ y \in X \mid y \leq x \}$ of $x$. 
A subset $A$ is a {\bf downset} {\rm(resp.} {\bf upset}) if $A = \bigcup_{x \in A} \mathop{\downarrow} x $ {\rm(resp.} $A = \bigcup_{x \in A} \mathop{\uparrow} x$). 
Define the height $\mathop{\mathrm{ht}} (x)$ of $x \in X$ by 
$
\mathop{\mathrm{ht}} (x) := \sup \{ |C| - 1 \mid C :\text{chain containing }x \text{ as the maximal point}\}
$. 
Define the height of the empty set is $-1$. 
The height $\mathop{\mathrm{ht}} (A)$ of a nonempty subset $A \subseteq X$ is defined by $\mathop{\mathrm{ht}} (A) := \sup_{x \in A} \mathop{\mathrm{ht}} (x)$.

\subsubsection{Abstract cell complex}
For a set $S$ with a transitive relation $\prec$ and a function $\dim \colon S \to \Z_{\geq 0}$, the triple $(S, \prec, \dim )$ is an {\bf abstract cell complex} if $x \prec y$ implies $\dim x < \dim y$. 
Then $\dim x$ is called the {\bf dimension} of $x$, and $x$ is call a {\bf cell}. 
A $k$-cell is a cell whose dimension is $k$. 
The {\bf codimension} of $x$ is $\sup_{y \in S} \dim y - \dim x$. 
For a finite pre-ordered set $(X, \leq)$, 
the triple $(X, <, \mathop{\mathrm{ht}})$ is an abstract cell complex. 

\subsection{Notion of topology}

A subset $A$ of a topological space is a {\bf collar} of a subset $B \subseteq A$ if there is a homeomorphism $h \colon [0,1] \times B \to A$ with $h(\{ 0\} \times B) = B$. 

A point $x$ of a topological space $X$ is $\bm{T_0}$ (or Kolmogorov) if for any point $y \neq x \in X$ there is an open subset $U$ of $X$ such that $|\{x, y \} \cap U| =1$, where $|A|$ is the cardinality of a subset $A$.
A topological space is $T_0$ if each point is $T_0$.

The {\bf specialization preorder} $\leq_{\tau}$ on a topological space $(X, \tau)$ is defined as follows: $ x \leq_{\tau} y $ if  $ x \in \overline{\{ y \}}$, where $\overline{\{ y \}}$ is the closure of $\{ y \}$. 
A topology $\tau$ is $T_0$ if and only if $\leq_{\tau}$ is a partial order. 
For a finite topological space, the subset is closed {\rm(resp.} open) if and only if it is a downset {\rm(resp.} upset) with respect to the specialization preorder. 
The height of a topological space, a point, and a subset are defined the height with respect to the specialization preorder. 
For any $k \in \Z_{\geq 0}$ and a topological space $X$, denote by $X_k$ the set of height $k$ points of $X$ and by $X_{\leq k}$ the set of points of $X$ whose height is less than or equal to $k$. 

\subsubsection{Surface graph}
By a {\bf surface}, we mean a paracompact two-dimensional manifold.
A {\bf graph} is a cell complex whose dimension is at most one and which is a geometric realization of an abstract multi-graph.
In other words, it can be drawn such that no edges cross each other.
Such a drawing is called a {\bf surface graph}.

\subsubsection{Double of a compact surface}

The {\bf double} $S_{\mathrm{dbl}}$ of a compact surface $S$ is the resulting closed surface by gluing two copies of $S$ as follows:  $S_{\mathrm{dbl}} := S \times \{+,-\} / \sim$ where $(x,+) \sim (x,-)$ if $x \in \partial S$.

\subsubsection{Curves and loops}
A {\bf curve} is a continuous mapping $C: I \to X$ where $I$ is a non-degenerate connected subset of a circle $\mathbb{S}^1$.
A curve is {\bf simple} if it is injective.
We also denote by $C$ the image of a curve $C$.
Denote by $\partial C := C(\partial I)$ the boundary of a curve $C$, where $\partial I$ is the boundary of $I \subset \mathbb{S}^1$. 
The set difference $\mathrm{int} C := C \setminus \partial C$ is called the interior of the curve $C$. 
A simple curve is a simple closed curve if its domain is $\mathbb{S}^1$ (i.e. $I = \mathbb{S}^1$).
A simple closed curve is also called a {\bf loop}. 
An {\bf arc} is a simple curve whose domain is an interval.

\subsection{Notion of dynamical systems}

A {\bf flow} on a topological space is a continuous $\mathbb{R}$-action. 
For any flow $v \colon  \R \times S \to S$, the subset $\{ v(t,x) \mid t  \in \R \}$ for any point $x \in S$ is called the orbit of $x$. 
A subset of $S$ is said to be {\bf invariant} (or {\bf saturated}) if it is a union of orbits. 
A flow $v$ on a topological space $S$ is {\bf topologically equivalent} to a flow $w$ on a topological space $T$ if there is a homeomorphism $h \colon S \to T$ such that the images of any orbits of $v$ are orbits of $w$ and that $h$ preserves the directions of orbits of $v$ and $w$. 

\subsubsection{Orbit arc and local topological equivalence}

A closed interval contained in an orbit is called a {\bf orbit arc}. 
The restriction of a flow $v$ to a subset $U$ of $S$ is {\bf topologically equivalent} to one of a flow $w$ to a subset $V$ on a topological space $T$ if there is a homeomorphism $h \colon U \to V$ such that the images of any orbit arcs of $v$ in $U$ are orbit arcs of $w$, the inverse images of any orbit arcs of $w$ in $V$ are orbit arcs of $v$, and $h$ and $h^{-1}$ preserve the directions of orbit arc of $v$ and $w$. 

\subsection{Fundamental concepts for vector fields}
Let $X$ be a vector field generating a flow $v_X$ on a surface $S$. 
The orbit of a point $x \in S$ for $v_X$ is denoted by $O(x)$ and called an orbit generated by $X$.
An orbit is {\bf singular} if it is a singleton, and is {\bf periodic} if it is a circle. 
An orbit is {\bf closed} if it is either singular or periodic. 
An orbit is non-singular if it is not singular.
Denote by $\bm{\mathop{\mathrm{Sing}}(X)}$ {\rm(resp.} $\bm{\mathop{\mathrm{Per}}(X)}$, $\bm{\mathop{\mathrm{Cl}}(X)}$) the union of singular {\rm(resp.} periodic, closed) orbits. 
A point $x \in S$ is {\bf singular} {\rm(resp.} non-singular, {\bf periodic}, {\bf closed}) if its orbit is singular {\rm(resp.} non-singular, periodic, closed). 
A singular point is {\bf isolated} if there is its \nbd which contains no singular points except it. 
Then the \nbd is called its isolated neighborhood. 

For a point $x \in S$, define the $\omega$-limit set $\omega(x)$ and the $\alpha$-limit set $\alpha(x)$ of $x$ as follows: $\omega(x) := \bigcap_{n\in \mathbb{R}}\overline{\{v_X(t,x) \mid t > n\}} $, $\alpha(x) := \bigcap_{n\in \mathbb{R}}\overline{\{v_X(t,x) \mid t < n\}} $, where $v_X \colon \R \times S \to S$ is the flow generated by $X$. 
A point $x \in S$ is {\bf recurrent} if $x \in \omega(x) \cup \alpha(x)$. 
Denote by $\bm{\mathrm{P}(X)}$ the set of non-recurrent points and by $\bm{\mathrm{R}(X)}$ the set of non-closed recurrent points. 
Then $S = \mathop{\mathrm{Cl}}(X) \sqcup \mathrm{P}(X) \sqcup \mathrm{R}(X)$. 

A vector field $X$ generating a flow $v_X$ is {\bf topologically equivalent} to a vector field $Y$ generating a flow $v_Y$ if $v_X$ is topologically equivalent to $v_Y$. 
The restriction to a subset $U$ of a vector field $X$ generating a flow $v_X$ is {\bf topologically equivalent} to the restriction to a subset $V$ of a vector field $Y$ generating a flow $v_Y$ if the restriction $v_X|_U$ is topologically equivalent to the restriction $v_Y|_V$.

\subsubsection{Transversals}

A $C^1$-curve $C$ on a suface $S$ is {\bf transverse} to a vector field $X$ if the tangent vector at $x \in C$ and $X(x)$ span the tangent space $T_X S$ (i.e. $T_x C + \R X(x) = T_x S$, where $\R X(x) = \{r X(x) \mid r \in \R \}$). 
Then $C$ is called a {\bf transversal} to $X$. 

\subsubsection{Indices of isolated singular points of vector fields on surfaces}

\begin{figure}
\begin{center}
\includegraphics[scale=0.15]{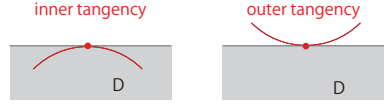}
\end{center}
\caption{Inner and outer tangencies.}
\label{fig:tangencies}
\end{figure}

A tangency $x$ of a curve $C$ to a loop bounding a closed disk $D$ is {\bf inner} (resp. {\bf outer}) if there is a small arc $I$ in $C$ containing $x$ such that the difference $I - \{ x \}$ is contained in the interior $\mathrm{int} D$ (resp. $(I - \{ x \}) \cap D = \emptyset$) as in left the (resp. right) on Figure~\ref{fig:tangencies}.
Recall the indices of isolated singular points of vector fields on surfaces (cf. \cite[\S~2.1.2]{2022VladislavMorse}). 
The {\bf index} of an isolated singular point outside of the boundary $\partial S$ is the difference $n_{i} - n_{o}$, where $n_i$ is the number of inner and $n_o$ is the number of outer tangencies of a loop which is transverse at all but finitely many points and bounds an open disk containing the singular point. 
Notice that the index of an isolated singular point outside of the boundary of a surface is independent of the choice of such a loop in the previous definition. 
For a singular point $x \in \partial S$, the {\bf index} $\bm{\mathrm{ind}_X(x)}$ of $X$ at $x$ is defined by $\mathrm{ind}_X(x) := \mathrm{ind}_{X_{\mathrm{dbl}}}(x)/2$, where $X_{\mathrm{dbl}}$ is the induced vector field on the double $S_{\mathrm{dbl}}$. 
We observe the following invariance of indices under small perturbations. 

\begin{lemma}\label{lem:inv_index}
Let $X$ be a vector field on a surface $S$ and $x \in S$ an isolated singular point. 
The following statements hold for any loop $C$ which is transverse at all but finitely many points and bounds an open disk $D$ whose closure is an isolated \nbd of $x$ and for any small positive number $\varepsilon >0$: 
\\
{\rm(1)} $\mathrm{ind}_X(x) = \sum_{y \in D \cap \mathop{\mathrm{Sing}}(Y)} \mathrm{ind}_Y(y)$ for any vector field $Y$ which is $C^0$-near $X$ and has only isolated singular points on $D$. 
\\
{\rm(2)} There is a loop $C' \subset B_{\varepsilon}(\partial D)$ which is transverse at all but finitely many points such that the number of the inner {\rm(resp.} outer{\rm)} tangencies on $C$ of $X$ equals one of the inner {\rm(resp.} outer) tangencies on $C'$ of $Y$, where $B_{\varepsilon}(\partial D)$ is the $\varepsilon$-\nbd of $\partial D$. 
\end{lemma}

\begin{proof}
By the invariance of transversality, any closed transverse arc $T$ for $X$ is also one for any vector field which is $C^0$-near to $X$. 

Fix a loop $C$ which is transverse at all but finitely many points and bounds an open disk $D$ whose closure is an isolated \nbd of an isolated singular point $x$. 
By arbitrarily small perturbation to $C$ near tangencies as in Figure~\ref{fig:perturbations_tang}, there are closed transverse arcs $T_j$ and a loop $\mu$ which is a finite disjoint union of orbit arcs $C_i$ and open subarcs $\gamma_j \subset T_j$. 
%
\begin{figure}
\begin{center}
\includegraphics[scale=0.6]{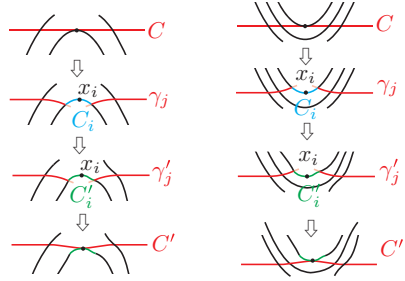}
\end{center}
\caption{Small perturbations near inner and outer tangencies.}
\label{fig:perturbations_tang}
\end{figure} 
Fix any points $x_i \in C_i$. 
Take a small $C^0$-\nbd $\mathcal{U}$ of $X$ such that for any vector field $Y$ there are orbit arcs $C'_i$ containing $x_i$ and connecting $\bigsqcup_{j} T_j$ and open subarc $\gamma'_j \subset T_j$ such that the disjoint union $\bigsqcup_{i} C'_i \sqcup \bigsqcup_{j} \gamma'_j$ is a loop $\mu'$ which is near $\mu$ with respect to the Hausdorff distance. 
Fix a small positive number $\varepsilon > 0$. 
Fix any $Y \in \mathcal{U}$. 
By arbitrarily small perturbation to $\mu'$ as in Figure~\ref{fig:perturbations_tang}, we obtain a loop $C'$ which is transverse at all but finitely many points with respect to $Y$ and bounds an open disk $D'$ such that $D$ and $D'$ are near with respect to the Hausdorff distance. 
Then $\partial D'$ is contained in the $\varepsilon$-\nbd $B_{\varepsilon}(\partial D)$ of $\partial D$. 
Therefore the symmetric difference $(D \cup D') - (D \cap D')$ contains no singular points of $X$ and $Y$. 
By construction, the loops $C = \partial D$ and $C' = \partial D'$ have the same numbers of inner and outer tangencies with respect to $X$ and $Y$ respectively. 
Since the Euler characteristic of a closed disk is one, Poincar{\'e}-Hopf theorem implies $\mathrm{ind}_X(x) = \sum_{y \in D \cap \mathop{\mathrm{Sing}}(X)} \mathrm{ind}_X(y) = \sum_{y \in D' \cap \mathop{\mathrm{Sing}}(Y)} \mathrm{ind}_Y(y) 
= \sum_{y \in D \cap \mathop{\mathrm{Sing}}(Y)} \mathrm{ind}_Y(y)$. 
\end{proof}

For any multi-saddle $x$ with an open disk $D$ as in the previous lemma, the intersection $D \cap \mathop{\mathrm{Sing}}(Y)$, which is a finite union of multi-saddles, is called  the {\bf continuation} $x_{Y}$ for $Y$ of $x$. 
Notice that the isolatedness of $Y$ in the previous lemma is not necessary if we define the index of an open disk whose boundary is a loop transverse at all but finitely many points.

\subsubsection{Centers}

An isolated singular point is a {\bf (topological) center} if there is its open \nbd to which the restriction of the flow is topologically equivalent to the flow generated by a vector field $X(x_1,x_2)=(-x_2,x_1)$ on the unit disk $\{ (x_1,x_2) \in \R^2 \mid x_1^2 + x_2^2 \leq 1 \}$. 
Denote by $\bm{\mathop{\mathrm{Sing}}_{\mathrm{c}}(X)}$ the set of centers of the flow generated by a vector field $X$.

\subsubsection{Multi-saddles}

\begin{figure}
\begin{center}
\includegraphics[scale=0.6]{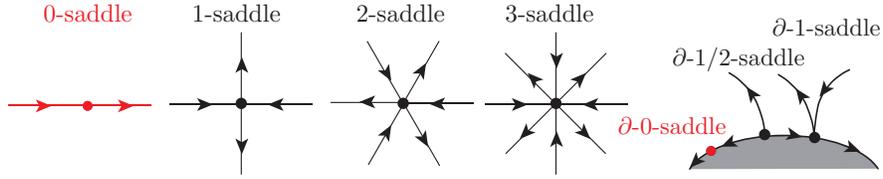}
\end{center}
\caption{Examples of multi-saddles}
\label{multi-saddles}
\end{figure}

A {\bf separatrix} is a non-singular orbit whose $\alpha$-limit or $\omega$-limit set is a singular point.
A {\bf $\bm{\partial}$-$\bm{k}$-saddle} {\rm(resp.} {\bf $\bm{k}$-saddle}) is an isolated singular point on {\rm(resp.} outside of) $\partial S$ with exactly $(2k + 2)$-separatrices, counted with multiplicity as in Figure~\ref{multi-saddles}.
A {\bf multi-saddle} is a $k$-saddle or a $\partial$-$(k/2)$-saddle for some $k \in \mathbb{Z}_{\geq 0}$.
A $1$-saddle is topologically an ordinary saddle and a $\partial$-$(1/2)$-saddle is topologically a $\partial$-saddle.
Notice that multi-saddles except $1$-saddles and $\partial$-$(1/2)$-saddles must be degenerate.
Here a singular point is degenerate if the Hessian (i.e. the determinant of the matrix of second partial derivatives at the point) is zero, and a singular point is nondegenerate if it is not degenerate. 

The {\bf multi-saddle connection diagram} $D(X)$ is the union of multi-saddles and separatrices from or to multi-saddles.
A {\bf multi-saddle connection} is a connected component of the multi-saddle connection diagram. 
It is known that any singular points of non-wandering flows with finitely many singular points on compact surfaces are either centers or multi-saddles \cite[Theorem 3]{cobo2010flows}. 
The indices of centers are one, and the indices of $k$-saddles and $\partial$-$k$-saddle are $-k$. 
Denote by $\bm{\mathrm{ind_+(X)}}$ the sum of indices of singular points which are positive for a vector field $X$ with finitely many singular points. 
Notice that the integer $\mathrm{ind_+(X)}$ for any non-wandering vector field $X$ with finitely many singular points on a compact surface corresponds to the number of centers.

\subsubsection{Self-connectedness}
A separatrix is {\bf self-connected} if either it is separatrix from and to a $k$-saddle or it connects between $\partial$-$k$-saddles on a boundary component. 
A multi-saddle connection is {\bf self-connected} if any separatrices contained in it are self-connected. 

\subsection{Fundamental concepts of Hamiltonian vector fields on a compact surface}

A $C^r$ vector field $Y$ for any $r \in \Z_{\geq 0} \sqcup \{ \infty \}$ on an orientable surface $\Sigma$ is {\bf Hamiltonian} if there is a $C^{r+1}$ function $H \colon \Sigma \to \mathbb{R}$ such that $dH= \omega(Y, \cdot )$ as a one-form, where $\omega$ is a volume form of $\Sigma$.
In other words, locally the Hamiltonian vector field $Y$ is defined by $Y = (- \partial H/ \partial x_2,  \partial H/ \partial x_1)$ for a local coordinate system $(x_1,x_2)$ of a point $p \in \Sigma$.
Note that a volume form on an orientable surface is a symplectic form.
It is known that local lipschitz continuity of vector fields guarantees the existence of their flows. 
On the other hand, Hamiltonian property of continuous vector fields also guarantees the existence of their flows. 
More precisely, one can observe the following statement, which is essentially mentioned in \cite{buczolich1994level} by applying a global implicit function theorem \cite[Theorem~1]{cater1990global}. 

\begin{lemma}
Any $C^0$ Hamiltonian vector field on a surface generates a flow. 
\end{lemma}

\begin{proof}
By Peano existence theorem~\cite{Peano1890diff}, the continuity of $X$ implies that the ordinal differential equation $dx/dt = X(x)$ for any $x \in \Sigma$ has a local solution. 
To complete this proof, it suffices to show the uniqueness of orbits. 
Fix a non-singular point $x_0 \in \Sigma$. 
Let $H$ be the $C^1$ Hamiltonian $H$ of $X$. 
Since the Hamiltonian $H$ of $X$ is $C^1$, there is a \nbd $U$ of $x_0$ such that the connected component of $U \cap H^{-1}(H(x_0))$ containing $x_0$ is the graph of a differentiable function, as mentioned in the introduction in \cite[p.27 lines~25 -- p.28 lines~4]{buczolich1994level}. 
Recall that $X(H) = 0$, because locally the Hamiltonian vector field $X$ is defined by $X = - \dfrac{\partial H}{\partial x_2} \dfrac{\partial}{\partial x_1} +  \dfrac{\partial H}{\partial x_1} \dfrac{\partial}{\partial x_2}$ for a local coordinate system $(x_1,x_2)$. 
This implies that the vector $X(x)$ for any non-singular point $x \in \Sigma$ is contained in the tangent space of the level set $H^{-1}(r)$ of $H$ for some $r \in \R$.  
Therefore the orbit of any non-singular point $x \in \Sigma$ is uniquely determined. 
This means that $X$ generates a flow. 
\end{proof}

It is known that a $C^r$ ($r \geq 1$) Hamiltonian vector field on a compact surface is structurally stable with respect to the set of $C^r$ Hamiltonian vector fields if and only if both each singular point is nondegenerate and each separatrix is self-connected (see  \cite[Theorem 2.3.8, p. 74]{ma2005geometric}).
Moreover, by \cite[Theorem 3]{cobo2010flows}, since the flow generated by any Hamiltonian vector field on a compact surface is non-wandering, each singular points of any Hamiltonian vector field with finitely many singular point on a compact surface is either a center or a multi-saddle such that a point is a center if and only if it is a local maximum or a local minimum of the Hamiltonian. 
%
In addition, one can not annihilate centers of any Hamiltonian vector field with finitely many singular points by any small perturbations with respect to the $C^0$-topology, because centers correspond to the local maximal values and local minimal values. 
%
Since every Hamiltonian vector field has no non-closed recurrent orbits, notice that the multi-saddle connection diagram of a Hamiltonian vector field on a compact surface is a complete invariant because the complement of the union of centers and the multi-saddle connection diagram is the disjoint union of open annuli consisting of periodic orbits (cf. \cite[Theorem~1.4.6]{ma2005geometric} and \cite[Corollary~6]{cobo2010flows}).
Moreover, any incompressible (i.e. divergence-free) vector fields on compact punctured spheres correspond to Hamiltonian vector fields up to topological equivalence (cf. \cite[Corollary C]{yokoyama2021relations}). 
%

%
%

%


\subsubsection{Spaces of Hamiltonian vector fields}

For any $r \in \Z_{\geq 0} \sqcup \{ \infty \}$ and any $g, p \in \mathbb{Z}_{\geq 0}$, let $\bm{\mathcal{H}^r(g,p)}$ be the set of $C^r$ Hamiltonian vector fields with finitely many singular points on a compact surface $\Sigma_{g,p}$ with genus $g$ and $p$ boundary components. 
For any distinct pair $(g,p) \neq (g',p')$, any flows generated by vector fields in $\mathcal{H}^r(g,p)$ and any flows generated by vector fields in $\mathcal{H}^r(g',p')$ belong to the different surfaces. 

Notice that $\mathcal{H}^r(g,p)$ can not capture the time evaluations of Hamiltonian flows with ``creations and annihilations of physical boundaries'' as in Figure~\ref{fig:creations_bdry} even if the number of centers (i.e. the sums of indices of singular points which are positive) is preserved, because ``creations and annihilations of physical boundaries'' as in Figure~\ref{fig:creations_bdry} change the base surfaces. 
To describe time evaluations of fluid phenomena in general settings, we need a topology which is different from the disjoint union topology on the set of $C^r$ Hamiltonian vector fields with finitely many singular points on compact surfaces. 
However, the construction of such a frame work is a future work. 
The paper will build the foundation for the framework.

By Gutierrez's smoothing theorem~\cite{gutierrez1978structural}, the non-existence of non-closed recurrent orbits of Hamiltonian vector fields implies the following smoothing properties.

\begin{lemma}\label{lem:smooth}
For any vector field $X$ in $\mathcal{H}^0(g,p)$, there is a vector field $Y$ in $\mathcal{H}^{\infty}(g,p)$ such that the flow generated by $X$ is topologically equivalent to the flow generated by $Y$. 
\end{lemma}

This implies that any pair of $\mathcal{H}^r(g,p)$ and $\mathcal{H}^s(g,p)$ for any $r, s \in \Z_{>0} \sqcup \{ \infty \}$ corresponds to each other up to topological equivalence. 
Put $\mathcal{H}(g,p) := \mathcal{H}^0(g,p)$.

\subsubsection{Creations and annihilations of singular points}


One can create new centers near degenerate singular points of any Hamiltonian vector field with degenerate singular points by arbitrarily small perturbations with respect to the $C^1$-topology. 
More precisely, for any Hamiltonian vector field $X$ with degenerate singular points in $\mathcal{H}^r(g,p)$, there is a Hamiltonian vector field $X' \in \mathcal{H}^r(g,p)$ which is topologically equivalent to $X$ such that for any \nbd $\mathcal{U} \subset \mathcal{H}^r(g,p)$ of $X'$ there is a Hamiltonian vector field $Y \in \mathcal{U}$ with $\mathrm{ind_+(X)} = \mathrm{ind_+(X')} < \mathrm{ind_+(Y)}$ by modifying a \nbd of the degenerate singular points.

To construct a foundation of description of time-evolving fluid phenomena,  we focus our analysis only on time evolutions that do not involve vortex generation and dissipation, such as vortex detachment in this paper.
Therefore we will exclude such instability of creations and annihilations of centers. 
More precisely, we restrict the space of Hamiltonian vector fields into the subspace by fixing the number of centers (i.e. the sums of indices of singular points which are positive) as follows. 

For any $r \in \Z_{\geq 0} \sqcup \{ \infty \}$ and any $i, g, p \in \mathbb{Z}_{\geq 0}$, let $\bm{\mathcal{H}^r(i,g,p)} \subseteq \mathcal{H}^r(g,p)$ be the set of $C^r$ Hamiltonian vector fields with finitely many singular points on the compact surface with genus $g$ and $p$ boundary components such that $\mathrm{ind_+(X)} = i$. 
%
%


\subsubsection{Fake multi-saddles}

A $\partial$-$k$-saddle {\rm(resp.} $k$-saddle) is a {\bf fake $\bm{\partial}$-saddle} {\rm(resp.} {\bf fake saddle}) if $k = 0$. 
A multi-saddle is {\bf fake} if it is either a $\partial$-$0$-saddle or a $0$-saddle. 
%
A singular point on the boundary is a {\bf pinching point}  \cite{sakajo2015transitions} if it is a $\partial$-$1$-saddle.
A Hamiltonian vector field with self-connected multi-saddle connections is {\bf $\bm{f}$-unstable} \cite{sakajo2015transitions} if it has just one fake multi-saddle and each singular point except the pinching point is non-degenerate. 
Here, $f$ in $f$-unstable represents ``fake''.

Generic perturbations of an $f$-unstable Hamiltonian vector field in $\mathcal{H}^r(i,g,p)$ imply a same structurally stable Hamiltonian vector field up to topological equivalence. 
More precisely, for any $f$-unstable Hamiltonian vector field $X \in \mathcal{H}^r(i,g,p)$, there is its \nbd $\mathcal{U} \subset \mathcal{H}^r(i,g,p)$ whose set difference removing $f$-unstable Hamiltonian vector fields consists of structurally stable Hamiltonian vector fields which are topologically equivalent to each other. 
Roughly speaking, a transition whose intermediate vector field is $f$-unstable is trivial under the non-existence of creations and annihilations of singular points. 

Therefore we will exclude such trivial instability of fake multi-saddles. 
More precisely, denote by $\bm{\mathcal{H}^r_*(i,g,p)} \subset \mathcal{H}^r(i,g,p)$ the set of Hamiltonian vector fields without fake multi-saddles on the compact surface with genus $g$ and $p$ boundary components such that $\mathrm{ind_+(X)} = i$. 
Then the subspace $\mathcal{H}^r_*(i,g,p)$ corresponds to the spaces without such trivial instability of fake multi-saddles. 
%
%
Moreover, we will show that the subspaces $\mathcal{H}_*^r(i,g,p)$ correspond to the spaces without creations and annihilations of singular points (see Lemma~\ref{lem:perturbation} for details). 

\subsubsection{Codimensions of multi-saddle connections and Hamiltonian vector fields}
Let $Y$ be a Hamiltonian vector field on a surface $T$ with finitely many singular points.  
The {\bf codimension} of a $k$-saddle $y$ ($k \geq 1$) is $2(k-1) = -2(1+ \mathrm{ind}(y))$. 
The {\bf codimension} of a $\partial$-$(l/2)$-saddle $y$ ($l \geq 1$) is $l-1 = -2(1/2 + \mathrm{ind}(y))$. 
Equivalently, the codimension of a $\partial$-$k$-saddle $y$ ($k \geq 1/2$) is $2k-1 = -2(1/2 + \mathrm{ind}(y))$. 
For a multi-saddle connection $D$ of $Y$, the codimension $\bm{\mathrm{codim}_{\mathrm{m}}(D)} \in \Z_{\geq 0}$ of $D$ with respect to multiplicity is defined as the sum of the codimensions of multi-saddles contained in $D$. 
The codimension $\bm{\mathrm{codim}_{\mathrm{h}}(D)} \in \Z_{\geq 0}$ of $D$ with respect to heteroclinicity is defined as the sum $N_{\partial}(D) + n_{m}(D) -1$, where $N_{\partial}(D)$ is the number of boundary components contained in $D$ and $n_{m}(D)$ is the number of multi-saddles outside of the boundary $\partial T$ contained in $D$. 

\begin{definition}
A multi-saddle connection $D$ has {\bf codimension $\bm{k}$} if $k = k_0 + \mathrm{codim}_{\mathrm{m}}(D) + \mathrm{codim}_{\mathrm{h}}(D)$, where $k_0$ is the number of fake multi-saddles contained in $D$. 
Then put $\bm{\mathrm{codim}(D)} := k$
\end{definition}
The multi-saddle connection diagram of $Y$ has {\bf codimension $\bm{k}$} if the sum of codimensions of multi-saddle connections is $k$. 
Then write $\bm{\mathrm{codim}(Y)} := k$. 
Moreover, denote by $\bm{\mathrm{codim}_{\mathrm{h}}(Y)}$ {\rm(resp.} $\bm{\mathrm{codim}_{\mathrm{m}}(Y)}$) the sum of the codimension with respect to heteroclinicity {\rm(resp.} multiplicity) of multi-saddle connections of $Y$.

\subsection{Codimensionality of subspaces of the space of Hamiltonian vector fields}

For any $k \in \Z_{\geq 0}$, denote by $\bm{\mathcal{H}^r_{*, k}(i,g,p)}$ the set of Hamiltonian vector field in $\mathcal{H}^r_*(i,g,p)$ whose multi-saddle connection diagram has codimension $k$. 
Note that codimension $k$ multi-saddle connections are also called atoms of complexity $k+1$ (cf. \cite[p. 94]{bolsinov2004integrable}). 
To make the hierarchical structure of instability of Hamiltonian vector fields correspond to the sum of local instability, we will use the codimension. 
In fact, the codimension of the multi-saddle connection diagram of a Hamiltonian vector field corresponds to the codimension of cells of the abstract cell complex structures of the space of Hamiltonian vector fields (see Theorem~\ref{main:02} for details). 

\subsubsection{High codimensional subspaces}

%
Write 
$\bm{\mathcal{H}^r_{*, >-1}(i,g,p)} := \mathcal{H}^r_*(i,g,p)$. 
For any $k \in \Z_{\geq 0}$, write $\bm{\mathcal{H}^r_{*, >k}(i,g,p)} := \bigsqcup_{l=k+1}^\infty \mathcal{H}^r_{*, l}(i,g,p)$. 
Then $\mathcal{H}^r_{*, >k}(i,g,p) = \mathcal{H}^r_{*, >k-1}(i,g,p) - \mathcal{H}^r_{*, k}(i,g,p)$ for any $k \in \Z_{\geq 0}$. 
From definition of $\mathcal{H}^r_*(i,g,p)$, we have the following relation: 
\[
\mathcal{H}^r_*(i,g,p) = \mathcal{H}^r_{*, >-1}(i,g,p) \supseteq \mathcal{H}^r_{*, >0}(i,g,p) \supseteq \mathcal{H}^r_{*, >1}(i,g,p) \supseteq \cdots
\]

\section{Perturbations}

\subsection{Description of the contiuations of the multi-saddle connections}

\subsubsection{Complete invariance of the union of the multi-saddle connection diagram and centers}

We observe the complete invariance of the union of the multi-saddle connection diagram and centers as follows.


\begin{lemma}\label{lem:fin_comb}
The union of the multi-saddle connection diagram and centers is a complete invariant of $\mathcal{H}^r(g,p)$ as a surface graph for any $r \in \Z_{\geq 0} \sqcup \{ \infty \}$ and any $g, p \in \mathbb{Z}_{\geq 0}$. 
Moreover, the complement of such a union is a finite disjoint union of invariant open periodic annuli. 
\end{lemma}

\begin{proof}
Let $X$ be a Hamiltonian vector field with finitely many singular points on the compact surface $\Sigma_{g,p}$ with genus $g$ and $p$ boundary components. 
From $\Sigma_{g,p} = \mathop{\mathrm{Cl}}(X) \sqcup \mathrm{P}(X) \sqcup \mathrm{R}(X)$, by the existence of the Hamiltonian of $X$, there are no non-closed recurrent orbits and so $\Sigma_{g,p} = \CX \sqcup \mathrm{P}(X)$. 
Since any Hamiltonian vector field is divergence-free and any divergence-free vector field on a compact manifold is non-wandering, the flow generated by $X$ is non-wandering. 
By \cite[Theorem 2.5 and Proposition 2.6]{yokoyama2019properness}, the union $\SX \sqcup \mathrm{P}(X) = \SX \cup D(X)$ is closed and the union $\PX$ is open dense. 
From \cite[Theorem 3]{cobo2010flows}, any singular points of a non-wandering flow with finitely many singular points on compact surfaces are either centers or multi-saddles. 
Then the union $\SX \cup D(X)$ of singular points and the multi-saddle connection diagram consists of finitely many orbits. 
Then $\Sigma_{g,p} - (\SX \cup D(X)) = \PX$. 
The openness of $\PX$ and the finiteness of $\SX \cup D(X)$ imply that the union $\PX$ is a finite disjoint union of invariant open periodic annuli. 
This means that the surface graph $\SX \cup D(X)$ is a complete invariant as a surface graph. 
\end{proof}

\subsubsection{Invariance of positive indexed singular points under small perturbations}

For Hamiltonian vector fields $X,Y$ with finitely many singular points on a compact surface $S$, the Hamiltonian vector field $Y$ {\bf has same type of positive indexed singular points} of $X$ if there are disjoint open subsets $U_1, \ldots , U_l$ of  $S$ satisfying the following conditions: 
\\
{\rm(1)} The disjoint union $\bigsqcup_{i=1}^l U_i$ is a \nbd of $\mathop{\mathrm{Sing}}(X) \cup \mathop{\mathrm{Sing}}(Y)$. 
\\
{\rm(2)} $|U_i \cap \mathop{\mathrm{Sing}}(X)| = 1 \leq |U_i \cap \mathop{\mathrm{Sing}}(Y)|$ for any $i \in \{1, \ldots , l \}$. 
\\
{\rm(3)} For any $i \in \{1, \ldots , l \}$, the singular point in $U_i \cap \mathop{\mathrm{Sing}}(X)$ is a center if and only if $|U_i \cap \mathop{\mathrm{Sing}}(Y)| = 1$ and the singular point in $U_i \cap \mathop{\mathrm{Sing}}(Y)$ is a center. 

%
%
%
%
We observe that any small perturbation can not merge a pair of distinct singular points for a vector field in $\mathcal{H}^r_*(i,g,p)$. 
More precisely, we have the following statement. 

\begin{lemma}\label{lem:perturbation}
For any $X \in \mathcal{H}^r_*(i,g,p)$, there is its \nbd in $\mathcal{H}^r_*(i,g,p)$ each element of which has same type of positive indexed singular points of $X$. 
\end{lemma}

\begin{proof}
Let $X$ be a Hamiltonian vector field with finitely many singular points on a compact surface $\Sigma_{g,p}$ with genus $g$ and $p$ boundary components and $H_X$ the Hamiltonian of $X$. 
By \cite[Theorem 3]{cobo2010flows}, any singular points of Hamiltonian vector fields with finitely many singular points on compact surfaces are either centers or multi-saddles. 
Then the set $\mathop{\mathrm{Sing}}_{\mathrm{c}}(X)$ of centers is the complement in $\mathop{\mathrm{Sing}}(X)$ of multi-saddles. 

Let $x_1, x_2, \ldots , x_{j}$ be the centers of $X$ and $x_{j+1}, x_{j+2}, \ldots , x_l$ the multi-saddles of $X$.  
Choose any small number $\varepsilon > 0$. 
Take open disks $V_1, V_2, \ldots , V_l$ which are isolated \nbds of $x_1, x_2, \ldots , x_{l}$ respectively and whose boundaries are loops transverse at all but finitely many points such that 
$\{ x \in \Sigma_{g,p} \mid \| X(x) \| < \varepsilon/2 \} \subset \bigsqcup_{j=1}^l V_j \subset \{ x \in \Sigma_{g,p} \mid \| X(x) \| < \varepsilon \}$. 
By Lemma~\ref{lem:inv_index}, there is a \nbd $\mathcal{U}$ of $X$ in $\mathcal{H}^r_*(i,g,p)$ such that $\mathrm{ind}_X(x_j) = \sum_{y \in V_j \cap \mathop{\mathrm{Sing}}(Y)} \mathrm{ind}_Y(y)$ and $\max_{x \in \Sigma_{g,p}} \| X(x) - Y(x) \| < \varepsilon /2$ for any $Y \in \mathcal{U}$. 
Then $V := \bigsqcup_{j=1}^l V_j$ is a \nbd of $\mathop{\mathrm{Sing}}(Y)$ for any $Y \in \mathcal{U}$. 

We claim that $V_1, V_2, \ldots , V_i$ are isolated \nbds of centers but that $V_{i+1}, V_{i+2}, \ldots , V_l$ contains no centers for any $Y \in \mathcal{U}$. 
Indeed, fix any $Y \in \mathcal{U}$. 
Since centers correspond to the local maximal values and local minimal values of the Hamiltonian $H_X$, any small perturbation for $X$ in $\mathcal{H}^r_*(i,g,p)$ does not decrease the number of centers. 
Therefore $V_1, V_2, \ldots , V_i$ contains at least one centers of $Y$. 
Since $Y \in \mathcal{H}^r_*(i,g,p)$, the vector field $Y$ contains at most $i$ centers. 
This means that $V_1, V_2, \ldots , V_i$ contain exactly one center of $Y$ and that $V_{i+1}, V_{i+2}, \ldots , V_l$ contain no centers of $Y$. 
\end{proof}

The previous lemma means that the subspaces $\mathcal{H}^r_*(i,g,p)$ correspond to the spaces without creations and annihilations of singular points because of Poincar{\'e}-Hopf theorem. 

\subsubsection{Existence of associated neighborhoods of the multi-saddle connections}

Write $\bm{N(X)}$ the number of multi-saddle connections of a vector field $X \in \mathcal{H}^r(g,p)$. 
Denote by $\bm{v_X(A)}$, called the {\bf saturation} of a subset $A$, the union of orbits of $X$ intersecting $A$.
We observe non-existence of merges of distinct multi-saddle connections as follows. 

\begin{figure}
\begin{center}
\includegraphics[scale=0.6]{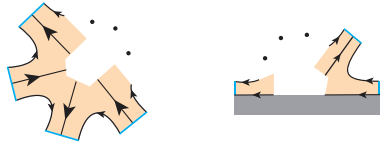}
\end{center}
\caption{A small collar of the boundary of a \nbd of multi-saddles outside of the boundary and one on the boundary.}
\label{multi-saddle_nbd}
\end{figure}

\begin{lemma}\label{lem:no_merge_msc}
For any $X \in \mathcal{H}^r(g,p)$, there are pairwise disjoint closed connected invariant \nbds $V_t$, called {\bf associated neighborhoods}, of the multi-saddle connections $D_{X,t}$
and a \nbd $\mathcal{U} \subset \mathcal{H}^r(g,p)$ of $X$ satisfying the following properties for any element $Y$ of $\mathcal{U}$: 
\\
{\rm(1)} The disjoint union $\bigsqcup_{t} V_t$ contains no periodic boundaries of $X$ but an open connected invariant \nbd $W(Y)$ of the multi-saddle connection diagram $D(Y)$ of $Y$ such that the set difference $W(Y) - D(Y)$ consists of finitely many open periodic annuli of $Y$ and that the multi-saddle connection diagram $D_{X,t}$ is a strong deformation retract of the intersection $W(Y) \cap V_t$. 
\\
{\rm(2)} Any connected components of the complement $\Sigma_{g,p} - \bigsqcup_{t} V_t$ contain open invariant periodic annuli of $Y$. 
\\
{\rm(3)} If $X \in \mathcal{H}^r_*(g,p)$, then $N(Y) \geq N(X)$ {\rm(i.e.} the number of multi-saddle connections of $X$ is more than or equals to one of $Y${\rm)}. 
\\
{\rm(4)} There are closed arcs $I_{t,j} \subset V_t$ which are transverse closed arcs and intersect separatrics of any $Y' \in \mathcal{U}$ such that any connected component $C_k$ of the set difference $W(Y') \setminus \bigsqcup_{t,j} I_{t,j}$ is a disk $B_k(Y')$ in which the restriction $Y'|_{A}$ of $Y'$ to a \nbd $A \subseteq B_k(Y)$ of $\partial C_k$ is topologically equivalent to the restriction to a small collar of the boundary of an isolated \nbd of either a $n$-saddle or $\partial$-$n/2$-saddle of the vector field $Z = (Z_x, Z_y) := ( r \cos (n \theta), r \cos (n \theta))$ as in Figure~\ref{multi-saddle_nbd}, where $(r,\theta)$ is the polar coordinate system. 
\end{lemma}

\begin{proof}
Let $U_1, \ldots , U_l$ be the connected components of the complement $\Sigma_{g,p} - (D(X) \sqcup \mathop{\mathrm{Sing}}_{\mathrm{c}}(X)) = \mathop{\mathrm{Per}}(X)$ and $T_1, \ldots , T_l$ transverse closed arcs in $U_1, \ldots , U_l$ respectively. 
Then $v_X(T_k)$ is a closed periodic annulus of $X$ in $U_k$ 
for any $k \in \{ 1, \ldots , l \}$. 
Fix total orders in $T_1, \ldots , T_l$ respectively. 

We claim that there are invariant closed annuli $U'_k \subset U_k$ and a \nbd $\mathcal{U} \subset \mathcal{H}^r(g,p)$ of $X$ such that for any $Y \in \mathcal{H}^r(g,p)$ any $U'_k$ contains an invariant open periodic annulus $A_k(Y)$ of $Y$ and is contained in an invariant open periodic annulus of $Y$ in $U_k$. 
Indeed, fix three points $x_{k,-} < x_k < x_{k,+} \in \mathop{\mathrm{int}} T_k$, pairwise disjoint closed  subarcs $T'_{k,-}, T'_k, T'_{k,+} \Subset T_k$ whose interiors contain $x_{k,-}$, $x_k$, $x_{k,+}$ respectively, where $A \Subset B$ means $\overline{A} \subseteq \mathop{\mathrm{int}} B$. 
The saturations $U'_{k,-} := v_X(T'_{k,-})$, $U'_k := v_X(T'_{k})$, and $U'_{k,+} := v_X(T'_{k,+})$ are pairwise disjoint invariant closed annuli in $U_k$ consisting of periodic orbits with respect to $X$. 
Denote by $T''_k \subset T_k$ the minimal closed interval which contains $T'_k \sqcup T'_{k,-} \sqcup T'_{k,+}$. 
The saturation $U''_k := v_X(T''_{k}) \subset U_k$ is a closed annulus consisting of periodic orbits with respect to $X$ such that the set diffrence $U''_k - (U'_{k,-} \sqcup U'_{k,+}) = v_X(T''_k - (T'_{k,-} \sqcup T'_{k,+}))$ is an invariant open annulus of $X$ containing $U'_k$. 
Choose a closed subarc $T'''_k \Subset T''_k$ whose interior contains $x_{k,-}$, $x_k$, $x_{k,+}$, and a closed subarc $T''''_k \Subset T'_k$ (resp. $T''''_{k,-} \Subset T'_{k,-}$, $T''''_{k,+} \Subset T'_{k, +}$) whose interior contains $x_k$ (resp. $x_{k,-}$, $x_{k,+}$). 
Then $T''''_k \Subset T'_k \Subset T'''_k \Subset T''_k \Subset T_k$. 
The saturations $U'''_k := v_X(T'''_{k}) \subset U_k$ and $U''''_k := v_X(T''''_{k}) \subset U_k$ are closed annuli consisting of periodic orbits with respect to $X$ such that $U''_k - (U'_{k,-} \sqcup U'_{k,+}) \subset U'''_k$. 
From the flow box theorem (cf. Theorem 1.1, p.45\cite{aranson1996introduction}), there are finitely many open flow boxes $B_{k,1}, \ldots , B_{k,k'} \subset U''_k$ whose union is a \nbd $V$ of the open periodic annulus $U'''_k$. 
Then there is a \nbd $\mathcal{U}_k$ of $X$ such that $T_k$ is a transversal for any element $Y \in \mathcal{U}_k$, that the positive orbit $O^+(x)$ of any point $x$ in $T'''_k$ (resp. $T''''_k$) for any element $Y \in \mathcal{U}_k$ intersects $T''_k$ (resp. $T'_k$) and that the orbit arcs from points in $T'''_k$ (resp. $T''''_k$, $T''''_{k,-}$, $T''''_{k,+}$) and to the first returns to $T_k$ are contained in $V \subset U''_k$ (resp. $U'_k$, $U'_{k,-}$, $U'_{k,+}$). 
Since any orbits of Hamiltonian vector fields are contained in level sets of the Hamiltonians and any distinct points in a transversal to a Hamiltonian vector field have different values of the Hamiltonian, the positive orbit of any point $x \in T'''_k$ for any element $Y \in \mathcal{U}_k$ intersects $T_k$ at itself. 
This means that any points in $T'''_k$ are periodic with respect to any element $Y \in \mathcal{U}_k$ such that $v_Y(T'''_k) \subset V \subset U''_k$, and that $v_Y(T''''_k) \subset U'_k$, $v_Y(T''''_{k,-}) \subset U'_{k,-}$ and $v_Y(T''''_{k,+}) \subset U'_{k,+}$. 
Therefore any $U'_k$ contains an invariant open periodic annulus $v_Y(T'''_k)$. 
Denote by $A_k(Y)$ the connected component of $U''_k - (v_Y(T''''_{k,-}) \sqcup v_Y(T''''_{k,+}))$ containing $x_k$. 
Then $A_k(Y) \subset U_k$ is an invariant open periodic annulus of $Y$. 
By $v_Y(T''''_{k,-}) \subset U'_{k,-} \subset U''_k$ and $v_Y(T''''_{k,+}) \subset U'_{k,+} \subset U''_k$, we have $x_k \in U'_k \subset U''_k - (U'_{k,-} \sqcup U'_{k,+}) \subset U''_k - (v_Y(T''''_{k,-}) \sqcup v_Y(T''''_{k,+}))$. 
Since $A_k(Y)$ is the connected component of $U''_k - (v_Y(T''''_{k,-}) \sqcup v_Y(T''''_{k,+}))$ containing $x$, the connected component $A_k(Y)$ contains $U'_k$. 
Then the intersection $\mathcal{U} := \bigcap_{k=1}^l \mathcal{U}_k$ and $U'_k$ are desired.  
This completes the claim. 

Fix such invariant closed periodic annuli $U'_k \subset U_k$ as in the previous claim.
Denote by $V_t$ the connected component of $\Sigma_{g,p} - \bigsqcup_{k=1}^l U'_k$ containing a multi-saddle connection $D_{X,t}$ of $X$. 
By construction, for any $Y \in \mathcal{U}$, the connected components $V_t$ intersect no periodic boundary components of $X$, and the disjoint union $\bigsqcup_{t} V_t$ contains an open connected invariant \nbd $W(Y) := \Sigma_{g,p} - \bigsqcup_{k=1}^l \overline{A_k(Y)}$ of the multi-saddle connection diagram $D(Y)$ of $Y$ such that the complement $\Sigma_{g,p} - (D(Y) \sqcup \bigsqcup_{k=1}^l \overline{A_k(Y)}) = W(Y) - D(Y)$ consists of finitely many open periodic annuli of $Y$. 
By construction, the multi-saddle connection diagram $D_{X,t}$ is a strong deformation retract of the intersection $W(Y) \cap V_t$. 
This completes the proofs of assertions {\rm(1)} and {\rm(2)}. 

Suppose that $X \in \mathcal{H}^r_*(g,p)$. 
By construction, the non-existence of fake multi-saddles implies that $N(Y) \geq N(X)$ for any $Y \in \mathcal{U}$. 
Therefore assertion {\rm(3)} holds. 

For any separatrics $O_{t,j} \subset V_t$ of $X$, take transverse closed arcs $I_{t,j} \subset V_t$ whose boundaries are contained in $\partial V_t$. 
By the invariance of transversality, taking $\mathcal{U}$ small, we may assume that any $I_{t,j}$ are also transverse closed arcs for any $Y \in \mathcal{U}$. 
Then any connected component $C_k$ of the set difference $W(Y) \setminus \bigsqcup_{t,j} I_{t,j}$ is a disk $B_k(Y)$ such that any small collar of $\partial C_k$ in $B_k(Y)$ is topologically equivalent to a small collar of the boundary of an isolated \nbd of a multi-saddle as in Figure~\ref{multi-saddle_nbd}. 
\end{proof}

The intersection $D(Y) \cap V_t$, which is a finite disjoint union of multi-saddle connections of $Y$, is called the {\bf continuation} $\bm{D_{Y,t}}$ for $Y$ of the multi-saddle connection $D_{X,t}$. 
Denote by $\bm{\mathrm{codim}(D_{Y,t})}$, called the codimension of $D_{Y,t}$, the sum of the codimensions of multi-saddle connections contained in the continuation $D_{Y,t}$ for $Y$ of $D_{X,t}$. 
Similarly, denote by $\bm{\mathrm{codim}_{\mathrm{m}}(D_{Y,t})}$ (resp. $\bm{\mathrm{codim}_{\mathrm{h}}(D_{Y,t})}$), called the codimension of $D_{Y,t}$ with respect to multiplicity (resp. heteroclinicity), the sum of the codimensions with respect to multiplicity (resp. heteroclinicity) of multi-saddle connections contained in the continuation $D_{Y,t}$ for $Y$ of $D_{X,t}$. 


\subsubsection{Hamiltonian vector fields with same type of multi-saddles}

From now on, we focus only $C^0$ perturbations in this section because small $C^r$ perturbations are small $C^0$ perturbations.  
Put $\mathcal{H}_*(i,g,p) := \mathcal{H}^0_*(i,g,p)$.

For Hamiltonian vector fields $X,Y$ with finitely many singular points on a compact surface $S$, the Hamiltonian vector field $Y$ {\bf has same type of multi-saddles} of $X$ if there are disjoint open subsets $U_1, \ldots , U_l$ of  $S$ satisfying the following conditions: 
\\
{\rm(1)} The disjoint union $\bigsqcup_{i=1}^l U_i$ is a \nbd of $\mathop{\mathrm{Sing}}(X) \cup \mathop{\mathrm{Sing}}(Y)$. 
\\
{\rm(2)} $|U_i \cap \mathop{\mathrm{Sing}}(X)| = 1$ for any $i \in \{1, \ldots , l \}$. 
\\
{\rm(3)} For any $i \in \{1, \ldots , l \}$, the singular point in $U_i \cap \mathop{\mathrm{Sing}}(X)$ is a $k$-saddle {\rm(resp.} $\partial$-$k/2$-saddle) if and only if $|U_i \cap \mathop{\mathrm{Sing}}(Y)| = 1$ and the singular point in $U_i \cap \mathop{\mathrm{Sing}}(Y)$ is a $k$-saddle {\rm(resp.} $\partial$-$k/2$-saddle). 

A boundary component (i.e. a connected component of the boundary) of a surface is {\bf periodic} if it is a periodic orbit. 
A Hamiltonian vector field $Y$ {\bf has same type of singular points} of a Hamiltonian vector field $X$ if $Y$ has same type of centers and same type of multi-saddles. 

For any $X \in \mathcal{H}_*(i,g,p)$, any multi-saddle connection $D_{X}$ of $X$ with an associated neighborhood $V_{D_{X}}$, any continuation $D_{Y}$ of $D_X$ for $Y \in \mathcal{H}_*(i,g,p)$, and a vector filed $Z$ which is either $X$ or $Y$, write the following symbols: 
\\
$\bm{n_{\partial}(D_{Z})}$ : the number of multi-saddles contained in $D_{Z}$ of $Z$ on the boundary $\partial \Sigma(g,p)$
\\
$\bm{n_{m}(D_{Z})}$ : the number of multi-saddles contained in $D_{Z}$ of $Z$ outside of $\partial \Sigma(g,p)$
\\
$\bm{n(D_{Z})}  := n_{\partial}(D_{Z}) + n_{m}(D_{Z})$ : the number of multi-saddles contained in $D_{Z}$ of $Z$  
 \\
$\bm{N(D_{Z})}$ : the number of multi-saddle connections contained in $D_{Z}$ of $Z$  
\\
$\bm{N_{\mathrm{per}}(D_{Z})}$ : the number of periodic boundary components contained in an associated neighborhood $V_{D_{X}}$
\\
$\bm{M(D_{Z,t})} := N_{\mathrm{per}}(D_{Z,t}) + N(D_{Z,t})$
\\
$\bm{N_{\partial}(D_{Z,t})}$ : the number of boundary components contained in $D_{Z,t}$
\\
We show the local maximality of codimensionality. 

\begin{lemma}\label{lem:deg_codim}
For any $X \in \mathcal{H}_*(i,g,p)$ and any multi-saddle connection $D_{X,t}$ of $X$, there is a \nbd $\mathcal{U} \subset \mathcal{H}_*(i,g,p)$ satisfying the following properties for any $Y \in \mathcal{U}$: 
\\
{\rm(1)} The following statements hold: 
\\
{\rm(1.1)} $\mathrm{codim}(D_{X,t}) - \mathrm{codim}(D_{Y,t}) = n_{\partial}(D_{Y,t}) - n_{\partial}(D_{X,t}) + n_{m}(D_{Y,t}) - n_{m}(D_{X,t}) + M(D_{Y,t}) - M(D_{X,t}) = n(D_{Y,t}) - n(D_{X,t}) + M(D_{Y,t}) - M(D_{X,t}) \geq 0$. 
\\
{\rm(1.2)} $n(D_{Y,t})  \geq n(D_{X,t})$, $n_{m}(D_{Y,t}) \geq n_{m}(D_{X,t})$, $N_{\mathrm{per}}(D_{Y,t}) \geq N_{\mathrm{per}}(D_{X,t})$, and $M(D_{Y,t}) \geq M(D_{X,t})$. 
\\
{\rm(2)} The following are equivalent: 

{\rm(a)} $n_{\partial}(D_{Y,t}) = n_{\partial}(D_{X,t})$ and  $n_{m}(Y) = n_{m}(D_{X,t})$ and $M(D_{Y,t}) = M(D_{X,t})$. 

{\rm(b)} $\mathrm{codim}(D_{Y,t}) = \mathrm{codim}(D_{X,t})$. 

{\rm(c)} The multi-saddle connections $D_{X,t}$ and $D_{Y,t}$ are isomorphic as a surface graph. 
%
\end{lemma}

\begin{proof}
By Lemma~\ref{lem:no_merge_msc} and Lemma~\ref{lem:perturbation}, take pairwise disjoint closed invariant small associated neighborhoods $V_t$ of multi-saddle connections $D_{X,t}$ intersecting no periodic boundary components and a \nbd $\mathcal{U} \subset \mathcal{H}_*(i, g,p)$ of the vector field $X$ as in Lemma~\ref{lem:no_merge_msc} such that each element of $\mathcal{U}$ has same type of positive indexed singular points of $X$. 
For any separatrices $O_{t,j} \subset V_t$ of $X$, take transverse closed arcs $I_{t,j} \subset V_t$ which intersects $O_{t,j}$ and whose boundaries are contained in $\partial V_t$ as in Lemma~\ref{lem:no_merge_msc}. 
Since the non-singular property on a compact subset is invariant under small perturbations, taking $\mathcal{U}$ small, we may assume that $N_{\mathrm{per}}(D_{X,t}) \leq N_{\mathrm{per}}(D_{Y,t})$ for any $Y \in \mathcal{U}$. 

Fix any $Y \in \mathcal{U}$. 
Since any connected components of $\bigsqcup_{t} V_t$ contains exactly one multi-saddle connection of $X$ and at least one multi-saddle connections of $Y$, we have that $N(D_{X,t}) = 1 \leq N(D_{Y,t})$ and so $M(D_{X,t}) \leq M(D_{Y,t})$ because $N_{\mathrm{per}}(D_{X}) = 0 \leq N_{\mathrm{per}}(D_{Y})$. 
From definition of codimension, the non-existence of fake multi-saddles implies that $\mathrm{codim}(D_{Z,t}) = \mathrm{codim}_{\mathrm{m}}(D_{Z,t}) + \mathrm{codim}_{\mathrm{h}}(D_{Z,t})$ for any $Z \in \{X,Y\}$, where 
\[
\begin{split}
\mathrm{codim}_{\mathrm{h}}(D_{Z,t}) &= \sum_{D' : \text{ multi-saddle connection in }D_{Z,t}} \left( N_{\partial}(D') + n_{m}(D') -1\right) 
\\ &= N_{\partial}(D_{Z,t}) + n_{m}(D_{Z,t}) - N(D_{Z,t})
\end{split}
\]
and 
\[
\mathrm{codim}_{\mathrm{m}}(D_{Z,t}) = \sum_{l} 2(l-1) |\{ l\text{-saddle in }D_{Z,t} \}| + \sum_{l} (l-1)|\{ \partial\text{-}(l/2)\text{-saddle in }D_{Z,t} \}|
\] 
is the sum of the codimensions of multi-saddles in $D_{Z,t}$. 
By $N_{\partial}(D_{X,t}) - N_{\partial}(D_{Y,t}) = N_{\mathrm{per}}(D_{Y,t}) - N_{\mathrm{per}}(D_{X,t}) \geq 0$, we obtain the following equation:  
\[
\begin{split}
&\mathrm{codim}_{\mathrm{h}}(D_{X,t}) - \mathrm{codim}_{\mathrm{h}}(D_{Y,t}) 
\\
&= N_{\partial}(D_{X,t}) - N_{\partial}(D_{Y,t}) + n_{m}(D_{X,t}) - n_{m}(D_{Y,t}) + N(D_{Y,t}) - N(D_{X,t})
\\
&= N_{\mathrm{per}}(D_{Y,t}) - N_{\mathrm{per}}(D_{X,t}) + n_{m}(D_{X,t}) - n_{m}(D_{Y,t}) + N(D_{Y,t}) - N(D_{X,t})
\\
&= M(D_{Y,t})) - M(D_{X,t})) + n_{m}(D_{X,t}) - n_{m}(D_{Y,t}) 
\end{split}
\]

We claim the following equalities: 
\[
\begin{split}
\mathrm{codim}_{\mathrm{m}}(D_{X,t}) - \mathrm{codim}_{\mathrm{m}}(D_{Y,t}) 
&= n_{\partial}(D_{Y,t}) -  n_{\partial}(D_{X,t}) + 2(n_{m}(D_{Y,t}) - n_{m}(D_{X,t}))
\\
&= n(D_{Y,t}) -  n(D_{X,t}) + n_{m}(D_{Y,t}) - n_{m}(D_{X,t})
\end{split}
\]
Indeed, for any $s,k \geq 1$ and for any $k$-saddle $x$ of $X$ and for any $k_j$-saddles $x_j$ of $Y$ with $k = \sum_{j=1}^s k_j$, we have that $\mathrm{codim}_X(x) = 2(k-1) \geq 2(k - s) = \sum_{j=1}^s 2(k_j -1) = \sum_{j=1}^s \mathrm{codim}_Y(x_j)$ and so $\mathrm{codim}_X(x) - \sum_{j=1}^s \mathrm{codim}_Y(x_j) = 2(s -1)$. 
For any $s', s'' \geq 0$ with $s' + s'' \geq 1$,  and any $k \geq 1$ and for any $\partial$-$k/2$-saddle $x$ of $X$ and for any $\partial$-$(k''_{j''}/2)$-saddles $x''_{j''}$ and any $k'_{j'}$-saddles $x'_{j'}$ of $Y$ with $k/2 = \sum_{j''=1}^{s''} k''_j/2 + \sum_{j'=1}^{s'} k'_{j'}$, we have that $\mathrm{codim}_X(x) = k-1 \geq k - s'' - 2s' = \sum_{j''=1}^{s''} (k_{j''} -1) + \sum_{j'=1}^{s'} 2(k'_{j'} -1) = \sum_{j''=1}^{s''} \mathrm{codim}_Y(x_{j''}) + \sum_{j'=1}^{s'} \mathrm{codim}_Y(x'_{j'})$ and so $\mathrm{codim}_X(x) - \sum_{j'=1}^{s'} \mathrm{codim}_Y(x'_{j'}) = s'' + 2s' -1$. 
Therefore we have the following inequality: 
\[
\begin{split}
0 \leq & \, \, \mathrm{codim}_{\mathrm{m}}(D_{X,t}) - \mathrm{codim}_{\mathrm{m}}(D_{Y,t}) 
\\
&= \sum_{t'=1}^{n_{m}(D_{X,t})} 2(s_{t'} -1) + \sum_{t''=1}^{n_{\partial}(D_{X,t})} (s''_{t''} + 2s'_{t''} -1)
\\
&= \sum_{t'=1}^{n_{m}(D_{X,t'})}  2s_{t'}  + \sum_{t''=1}^{n_{\partial}(D_{X,t})} (s''_{t''} + 2s'_{t''}) -  (2 n_{m}(D_{X,t}) + n_{\partial}(D_{X,t}))
\\
& = 2 n_{m}(D_{Y,t}) + n_{\partial}(D_{Y,t}) -  (2 n_{m}(D_{X,t}) + n_{\partial}(D_{X,t}))
\\
& = n_{\partial}(D_{Y,t}) - n_{\partial}(D_{X,t}) + 2(n_{m}(D_{Y,t}) - n_{m}(D_{X,t}))
\\
&= n(D_{Y,t}) -  n(D_{X,t}) + n_{m}(D_{Y,t}) - n_{m}(D_{X,t})
\end{split}
\]
because $\sum_{t''=1}^{n_{\partial}(D_{X,t})} s''_{t''}$ {\rm(resp.} $\sum_{t'=1}^{n_{m}(D_{X,t'})} s_{t'}+ \sum_{t''=1}^{n_{\partial}(D_{X,t})} s'_{t''}$) is the number of multi-saddles of $Y$ on {\rm(resp.} outside of) $\partial \Sigma(g,p)$ in $V_t$. 

By construction, we obtain that $n(D_{Y,t}) \geq  n(D_{X,t})$, $N(D_{Y,t}) \geq  N(D_{X,t})$ and $n_{m}(D_{Y,t}) \geq n_{m}(D_{X,t})$. 
This implies the assertion {\rm(1)}.

Fix an open invariant \nbd $W(Y) \subset \bigsqcup_{t} V_t$ of the multi-saddle connection diagram $D(Y)$ of $Y$ as in Lemma~\ref{lem:no_merge_msc}. 

If the multi-saddle connections $D_{X,t}$ and $D_{Y,t}$ are isomorphic as a surface graph, then 
the assertions {\rm(a)} and  {\rm(b)} hold. 

Suppose that $n_{\partial}(Y) = n_{\partial}(Y)$ and $n_{m}(Y) = n_{m}(X)$ and $M(Y) = M(X)$. 
By $M(Y) = M(X)$, each $V_t$ contains exactly one multi-saddle connection $D_{Y,t}$ of $Y$. 

We claim that we may assume that $Y$ has the same type of singular points of $X$. 
Indeed, let $x_1, \ldots , x_r$ be the multi-saddles of $X$. 
By Lemma~\ref{lem:inv_index}, by taking $\mathcal{U}$ small if necessary, there are pairwise disjoint isolated \nbds $D(x_j)$ of $x_j$ such that $\bigsqcup_{j=1}^r D(x_j)$ is a \nbd of the multi-saddles of any element of $\mathcal{U}$. 
From $n_{\partial}(Y) = n_{\partial}(Y)$ and $n_{m}(Y) = n_{m}(X)$, any $D(x_j)$ contains exactly one singular point $x'_j$ which is a multi-saddle such that $x_j$ is on $\partial \Sigma(g,p)$ if and only if $x'_j$ is on $\partial \Sigma(g,p)$. 
By Lemma~\ref{lem:inv_index}, we obtain $\mathrm{ind}_X(x_j) = \sum_{y \in D(x_j) \cap \mathop{\mathrm{Sing}}(Y)} \mathrm{ind}_Y(y) = \mathrm{ind}_Y(x'_j)$. 
Therefore $x_j$ and $x'_j$ are multi-saddles of the same type.
This means that $Y$ has same type of singular points of $X$. 

We claim that $D_{Y,t}$ is isomorphic to $D_{X,t}$ as a surface graph for any $i$. 
Indeed, fix a multi-saddle connections $D_{Y,t} \subset V_t$. 
Then it suffices to show that $D_{Y,t}$ is isomorphic to $D_{X,t}$ as a surface graph. 
Notice that any multi-saddles in $D_{Y,t}$ have a same values of the Hamiltonian of $Y$. 
For any connected component $B_{k}(Y)$ of the set difference $W(Y) \setminus \bigsqcup_{t,j} I_{t,j}$,  the restriction $Y|_{B_{k}(Y)}$ is topologically equivalent to the restriction $X|_{B_k(X)}$, where $B_k(X)$ is the connected component of $V_t \setminus \bigsqcup_{t,j} I_{t,j}$ containing the unique multi-saddle of $X$ in $B_{k}(Y)$ as in Figure~\ref{multi-saddle_nbd_std}. 
\begin{figure}
\begin{center}
\includegraphics[scale=0.6]{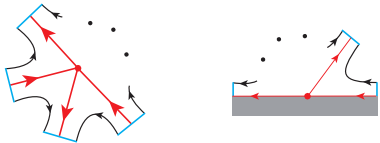}
\end{center}
\caption{A small collar of the boundary of an isolated \nbd of a multi-saddle outside of the boundary and one on the boundary.}
\label{multi-saddle_nbd_std}
\end{figure} 
Since any distinct points in a transverse closed arc have different values of the Hamiltonian of $Y$, the surface graph $D_{Y,t}$ is obtained by gluing the disks $B_{k}(Y)$ along $I_{t,j}$ as same as $D_{X,t}$. 
This means that $D(Y)$ is isomorphic to $D_{X,t}$ as a surface graph. 

The previous claim implies that the multi-saddle connections $D_{X,t}$ and $D_{Y,t}$ are isomorphic as a surface graph. 


Suppose that $\mathrm{codim}(D_{Y,t}) = \mathrm{codim}(D_{X,t})$. 
Since the differences $n(D_{Y,t}) - n(D_{X,t})$ and $M(D_{Y,t}) - M(D_{X,t})$ are non-negative, we have $n(D_{Y,t}) = n(D_{X,t})$ and $M(D_{Y,t}) = M(D_{X,t})$. 
By $N(D_{Y,t}) \geq N(D_{X,t})$ and $N_{\mathrm{per}}(D_{Y,t}) \geq N_{\mathrm{per}}(D_{X,t})$, the definition of $M(D_{Y,t}) = N(D_{Y,t}) + N_{\mathrm{per}}(D_{Y,t})$ implies that $N(D_{Y,t}) = N(D_{X,t}) = 1$ and $N_{\mathrm{per}}(D_{Y,t}) = N_{\mathrm{per}}(D_{X,t}) = 0$. 
Then $N_{\partial}(D_{X,t}) - N_{\partial}(D_{Y,t}) = N_{\mathrm{per}}(D_{Y,t}) - N_{\mathrm{per}}(D_{X,t}) = 0$ and $D_{Y,t}$ is a multi-saddle connection of $Y$. 
The non-existence of fake multi-saddles implies that the indices of any singular points are nonzero. 
Lemma~\ref{lem:inv_index} implies that $s_t = 1$ or either $\{ s'_t, s''_t \} = \{ 0, 1\}$ and so that the number of separatrices in $D_{X,t}$ equal one in $D_{Y,t}$. 
Lemma~\ref{lem:no_merge_msc} implies that the transverse closed arcs $I_{t,j}$ intersect exactly one separatrix of $D_{Y,t}$ as same as $D_{X,t}$. 
We claim that $s'_t \neq 1$. 
Indeed, assume that $s'_t = 1$. 
There is a $k$-saddle $x_Y \in D_{Y,t}$ of $Y$ which is a continuation of a $\partial$-$k$-saddle $x$ of $X$ on a boundary component $\partial_x$. 
By construction, we have that $x_Y \notin \partial_x$ and no separatrices of $x_Y$ with resect to $Y$ intersect $\partial_x$. 
From $N_{\mathrm{per}}(D_{Y,t}) = 0$, there are no new periodic boundary components on $Y$ in $V_t$ and so $\partial_x \subset D_{Y,t}$. 
Then there is a transverse closed arc $I_{t,j}$ which intersects a separatrix of $x$ in $\partial_x$ with respect to $X$. 
By Lemma~\ref{lem:no_merge_msc}, there is a separatrix of $x_Y$ of $Y$ which intersects $I_{t,j}$. 
Since no separatrices of $x_Y$ with resect to $Y$ intersect $\partial_x$, we obtain that $x_Y$ and $\partial_x$ intersect different points in $I_{t,j}$. 
This implies that $x_Y$ and $\partial_x$ have different values of a Hamiltonian of $Y$, which contradicts that $x_Y$ and $\partial_x$ are contained in the multi-saddle connection $D_{Y,t}$ of $Y$. 
Thus $s_t = 1$ or $s''_t = 1$. 
Then the surface graph $D_{Y,t}$ is obtained by gluing the disks $B_{k}(Y)$ along $I_{t,j}$ as same as $D_{X,t}$ because any distinct points in a transverse closed arc have different values of the Hamiltonian of $Y$. 
This means that $D_{X,t}$ and $D_{Y,t}$ are isomorphic as a surface graph.  
%
\end{proof}

We use the following symbols for any Hamiltonian vector field $Z \in \mathcal{H}_*(i,g,p)$: 
\\
$\bm{n_{\partial}(Z)}$ : the number of multi-saddles of $Z$ on the boundary $\partial \Sigma(g,p)$
\\
$\bm{n_{m}(Z)}$ : the number of multi-saddles of $Z$ outside of $\partial \Sigma(g,p)$
\\
$\bm{n(Z)} := n_{\partial}(Z) + n_{m}(Z)$ : the number of multi-saddles of $Z$
\\
$\bm{N(Z)}$ : the number of multi-saddle connections of $Z$
 \\
 $\bm{N_{\mathrm{per}}(Z)}$ : the number of periodic boundary components of $Z$
 \\
 $\bm{M(Z)} := N_{\mathrm{per}}(Z) + N(Z)$
\\
$\bm{N_{\partial}(Z)} := p - N_{\mathrm{per}}(Z)$ : the number of boundary components contained in $D_{Z,t}$

Lemma~\ref{lem:fin_comb} and the previous lemma imply the following statements. 

\begin{corollary}\label{cor:deg_codim}
For any Hamiltonian vector field $X \in \mathcal{H}_*(i,g,p)$, there is a \nbd $\mathcal{U} \subset \mathcal{H}_*(i,g,p)$ satisfying the following properties for any $Y \in \mathcal{U}$: 
\\
{\rm(1)} The following statements hold: 
\\
{\rm(1.1)} $\mathrm{codim}(X) - \mathrm{codim}(Y) = n_{\partial}(Y) - n_{\partial}(X) + n_{m}(Y) - n_{m}(X) + M(Y) - M(X) = n(Y) - n(X) + M(Y) - M(X) \geq 0$. 
\\
{\rm(1.2)} $N(Y)  \geq N(X)$, $n(Y)  \geq n(X)$, $n_{m}(Y) \geq n_{m}(X)$, and $N_{\mathrm{per}}(Y) \geq N_{\mathrm{per}}(X)$, and $M(Y) \geq M(X)$.
\\
{\rm(2)} The following are equivalent: 

{\rm(a)} $n_{\partial}(Y) = n_{\partial}(X)$ and  $n_{m}(Y) = n_{m}(X)$ and $M(Y) = M(X)$. 

{\rm(b)} $\mathrm{codim}(Y) = \mathrm{codim}(X)$. 

{\rm(c)} $X$ and $Y$ are topologically equivalent. 
%
\end{corollary}

\begin{proposition}\label{prop:perturbation}
For any Hamiltonian vector field $X \in \mathcal{H}_*(i,g,p)$, there is its \nbd $\mathcal{U}$ in $\mathcal{H}_*(i,g,p)$ satisfying the following properties for any $Y \in \mathcal{U}$:
\\
{\rm(1)} $\mathrm{codim}(X) \geq \mathrm{codim}(Y)$. 
\\
{\rm(2)} $\mathrm{codim}(X) = \mathrm{codim}(Y)$ if and only if $X$ and $Y$ are topologically equivalent. 
\end{proposition}

The following proposition implies the following statement.

\begin{corollary}\label{lem:open}
For any $r \in \Z_{\geq 0} \sqcup \{ \infty \}$ and any $k, i,g,p \in \Z_{\geq 0}$, the subset $\mathcal{H}^r_{*,k}(i,g,p)$ is open in the subspace $\mathcal{H}^r_{*, >k-1}(i,g,p)$ and any vector field in the subspace $\mathcal{H}^r_{*,k}(i,g,p)$ is $C^r$-structurally stable in the subspace $\mathcal{H}^r_{*, >k-1}(i,g,p)$. 
\end{corollary}

\begin{proof}
Since any $C^0$-\nbd contains a $C^r$-neighborhood and $\mathcal{H}^r_{*, 0}(i,g,p) \subseteq \mathcal{H}^0_*(i,g,p)$, the openness and $C^r$-structural stability follow from Proposition~\ref{prop:perturbation}.
%
\end{proof}

\subsection{Existence of Hamiltonian vector fields with one codimension smaller}
\subsubsection{Hamiltonian vector fields obtained by a Whitehead move}


A collapsing of a heteroclinic separatrix and the inverse operation as in Figure~\ref{splitting} are called the {\bf Whitehead moves}. 
We have the following descriptions of Whitehead moves. 

\begin{lemma}\label{lem:perturbation_03}
For any Hamiltonian vector field $Y \in \mathcal{H}_*(i,g,p)$ containing a $k$-saddle {\rm(}$k>1${\rm)} with its isolated \nbd $U$ and for any integer $l \in \{ 1,2, \ldots , k-1 \}$, there is a Hamiltonian vector field $Y' \in \mathcal{H}_*(i,g,p)$ obtained by its arbitrarily small perturbation in $\mathcal{H}_*(i,g,p)$ with $\mathrm{codim}(Y) - \mathrm{codim}(Y') = 1$ such that the closure $\overline{U}$ contains exactly two $(k-l)$-saddle and an $l$-saddle with a non-self-connected separatrix between the pair as in the top on Figure~\ref{splitting}.


Similarly, for any Hamiltonian vector field $X \in \mathcal{H}_*(i,g,p)$ containing a $\partial$-$k/2$-saddle {\rm(}$k>1${\rm)} with its isolated \nbd $U$ and for any integer $l \in \{ 1,2, \ldots , k-1 \}$, there is a vector field $X' \in \mathcal{H}_*(i,g,p)$ {\rm(resp.} $X'' \in \mathcal{H}_*(i,g,p)${\rm)} obtained by its arbitrarily small perturbation in $\mathcal{H}_*(i,g,p)$ with $\mathrm{codim}(X) - \mathrm{codim}(X') = 1$ {\rm(resp.} $\mathrm{codim}(X) - \mathrm{codim}(X'') = 1${\rm)} such that the closure $\overline{U}$ contains exactly two singular points of $X'$ {\rm(resp.} $X''${\rm)}, which are a $\partial$-$(k-l)/2$-saddle and a multi-saddle of $X'$ {\rm(resp.} $X''${\rm)} whose index is $l/2$, and contains a heteroclinic separatrix between the pair as in the middle {\rm(resp.} bottom{\rm)} on Figure~\ref{splitting}. 
%
\end{lemma}

\begin{proof}
Let $H_X$ {\rm(resp.} $H_Y$) be the Hamiltonian of $X$ {\rm(resp.} $Y$).
Considering a small isolating \nbd of a $\partial$-$k/2$-saddle {\rm(resp.} $k$-saddle) and replacing the restriction of the Hamiltonian $H_X$ {\rm(resp.} $H_Y$) to the neighborhood, we can replace a $\partial$-$k/2$-saddle {\rm(resp.} $k$-saddle) into a pair of a $\partial$-$(k-l)/2$-saddle and a multi-saddle whose index is $l/2$ {\rm(resp.} a $(k-l)$-saddle and an $l$-saddle) with a heteroclinic {\rm(resp.} non-self-connected) separatrix between the pair as in Figure~\ref{splitting} while the resulting Hamiltonian is $C^1$-near $H_X$ {\rm(resp.} $H_Y$).
Therefore the resulting Hamiltonian vector field is $C^0$-near $X$ {\rm(resp.} $Y$).
\begin{figure}
\begin{center}
\includegraphics[scale=0.6]{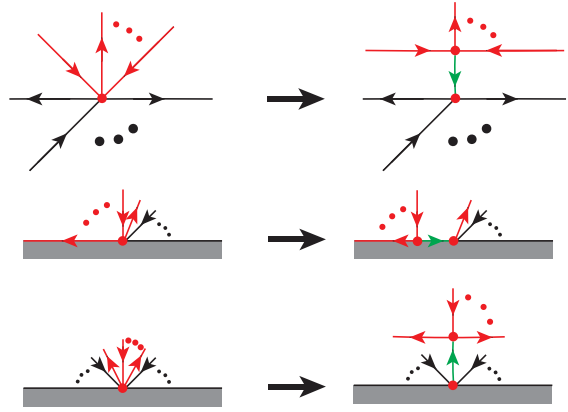}
\end{center}
\caption{Whitehead moves: splittings of multi-saddles}
\label{splitting}
\end{figure} 
\end{proof}

We call $X', X'', Y'$ in the previous lemma the Hamiltonian vector fields {\bf obtained by a Whitehead move}. 

\subsubsection{Hamiltonian vector fields obtained by pinching of multi-saddles on a same boundary}

For Hamiltonian vector fields $X,Y$ with finitely many singular points on a compact surface $S$ such that $\mathop{\mathrm{codim}}(X) - \mathop{\mathrm{codim}}(Y) =1$, the Hamiltonian vector field $Y$ is {\bf obtained by pinching of multi-saddles on a same boundary} if there are multi-saddles $\partial$-$k_0/2$-saddle $x_0$, $\ldots$ , $\partial$-$k_{l'}/2$-saddle $x_l'$ on a same boundary of $X$ and there are disjoint open subsets $U_0, U_1, \ldots , U_l$ of $S$ with $x_0 \in U_0$, $\ldots$ , $x_{l'} \in U_{l'}$ satisfying the following conditions: 
\\
{\rm(1)} The disjoint union $\bigsqcup_{i=0}^l U_i$ is a \nbd of $\mathop{\mathrm{Sing}}(X) \cup \mathop{\mathrm{Sing}}(Y)$. 
\\
{\rm(2)} $|U_i \cap \mathop{\mathrm{Sing}}(X)| = |U_i \cap \mathop{\mathrm{Sing}}(Y)| = 1$ for any $i \in \{l'+1, \ldots , l \}$. 
\\
{\rm(3)} For $i \in \{0, \ldots , l' \}$, the singular point in $U_{i} \cap \mathop{\mathrm{Sing}}(Y)$ is a $k_{i}/2$-saddle. 
\\
{\rm(4)} For any $i \in \{l'+1, \ldots , l \}$, the singular point in $U_i \cap \mathop{\mathrm{Sing}}(X)$ is a $k$-saddle {\rm(resp.} $\partial$-$k/2$-saddle) if and only if $|U_i \cap \mathop{\mathrm{Sing}}(Y)| = 1$ and the singular point in $U_i \cap \mathop{\mathrm{Sing}}(Y)$ is a $k$-saddle {\rm(resp.} $\partial$-$k/2$-saddle).

We have the following description of splittings of multi-saddle connections. 

\begin{figure}
\begin{center}
\includegraphics[scale=0.2]{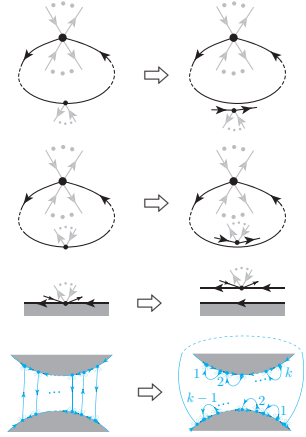}
\end{center}
\caption{Splittings of multi-saddle connections}
\label{splitting_02}
\end{figure}

\begin{lemma}\label{lem:perturbation_04}
Let $X \in \mathcal{H}^r_*(i,g,p)$ be a Hamiltonian vector field for any integer $r \in \Z_{\geq 0} \sqcup \{ \infty \}$. 
The following statements hold for a multi-saddle connection $D$ of $X$: 
\\
{\rm(1)} If $r \neq \infty$ and the multi-saddle connection $D$ contains at least two multi-saddles at least one of which is outside of the boudnary, then there is a Hamiltonian vector field $X' \in \mathcal{H}^r_*(i,g,p)$ obtained by its arbitrarily small $C^r$-perturbation in $\mathcal{H}^r_*(i,g,p)$ with $\mathrm{codim}(X) - \mathrm{codim}(X') = 1 = N(X') - N(X) = M(X') - M(X)$ such that $X'$ has same type of singular points of $X$ as in the top two on Figure~\ref{splitting_02}. 
\\
{\rm(2)} If $r = 0$ and the multi-saddle connection $D$ contains a $\partial$-$k$-saddle {\rm(}$k \in \Z_{\geq 1}${\rm)}, then there is a vector field $X'' \in \mathcal{H}_*(i,g,p)$ obtained by its arbitrarily small perturbation in $\mathcal{H}^0_*(i,g,p)$ with $\mathrm{codim}(X) - \mathrm{codim}(X'') = 1 = M(X'') - M(X)$ as in the second from bottom on Figure~\ref{splitting_02} such that the support of $X - X''$ is an isolated \nbd of the $\partial$-$k$-saddle of $X$. 
\\
{\rm(3)} If the multi-saddle connection $D$ contains at least two boundary components, then there is a vector field $X''' \in \mathcal{H}^r_*(i,g,p)$ obtained by its arbitrarily small $C^r$-perturbation in $\mathcal{H}^r_*(i,g,p)$ with $\mathrm{codim}(X) - \mathrm{codim}(X''') = 1 = N(X''') - N(X) = M(X''') - M(X)$ such that $X'''$ has same type of singular points of $X$ as in the bottom on Figure~\ref{splitting_02}. 
\end{lemma}

\begin{proof}
Suppose that the multi-saddle connection $D$ contains at least two multi-saddles $x, x'$ with $x \notin \partial \Sigma(g,p)$ (resp. a $\partial$-$k$-saddle $x$ {\rm(}$k \in \Z_{\geq 1}${\rm)}). 
Considering a small isolating \nbd of $x$ and and replacing the restriction of the Hamiltonian $H_X$ to the neighborhood as in the top two (resp. second from bottom) on Figure~\ref{splitting_02}, we can replace $D$ into two multi-saddle connections (resp. either two multi-saddle connections or a pair of a multi-saddle connection and a periodic orbit on the boundary) such that the replacement is  arbitrarily $C^{r+1}$ (resp. $C^1$) small perturbation of the Hamiltonian which changes only a small \nbd of $x$.
Therefore the resulting Hamiltonian vector field $X'$ (resp. $X''$) is $C^{r}$-near (resp. $C^0$-near) $X$.
By construction, the vector field $X'$ (resp. $X''$) is desired. 

Suppose that the multi-saddle connection $D$ contains at least two boundary components $\partial_1$ and $\partial_2$. 
Put $B_{s}(\partial_1) := \{ x \in \Sigma(g,p) \mid d(x,\partial_1) \leq s \}$, where $d(x,\partial_1) := \min_{y \in \partial_1} d(x,y)$. 
Choose a small number $\varepsilon > 0$ such that $B_{\varepsilon}(\partial_1) \cap \partial_2 = \emptyset$ and that $\delta := \min \{ |X(y)| \mid y \in B_{\varepsilon}(\partial_1) - B_{\varepsilon/2}(\partial_1) \} > 0$. 
For a positive number $s$, define a smooth decreasing function $h \colon \R_{\geq 0} \to [0,1]$ with $h|_{[0,\varepsilon/2]} = 1$ and $h|_{[\varepsilon, \infty)} = 0$ and a smooth Hamiltonian $H_s \colon \Sigma(g,p) \to \R_{\geq 0}$ by $H_s(x) := s \cdot h(d(x,\partial_1))$. 
Denote by $X_{H_s}$ the Hamiltonian vector field of $H_s$. 
Then there is a small number $\varepsilon' > 0$ such that $\max_{x\in \Sigma(g,p)}|X_{H_{\varepsilon'}}(x)| = \max_{x\in B_{\varepsilon}(\partial_1) - B_{\varepsilon/2}(\partial_1)}|X_{H_{\varepsilon'}}(x)| < \delta$. 
Therefore the resulting Hamiltonian vector field $X''' := X + X_{H_{\varepsilon'}}$, whose Hamiltonian is $H_X + H_{\varepsilon'}$, has the same type of singular points of $X$ and is $C^r$-near $X$ and as in the bottom on Figure~\ref{splitting_02} becuase $H_{X'''}(\partial_2) = H_X(\partial_2) = H_X(\partial_1) \neq H_{X'''}(\partial_1)$.
By construction, the vector field $X'''$ is desired. 
\end{proof}

Notice that we can deform the replacements in the previous proof into a smooth one up to topologically equivalent. 
Therefore we have the following observation. 

\begin{lemma}\label{lem:perturbation_04+}
Let $X_0 \in \mathcal{H}^r_*(i,g,p)$ be a Hamiltonian vector field for any integer $r \in \Z_{\geq 0}  \sqcup \{ \infty \}$ and $D$ a multi-saddle connection of $X$. 
There is a smooth Hamiltonian vector field $X \in \mathcal{H}^r_*(i,g,p)$ which is topologically equivalent to $X_0$ satisfying the following statements: 
\\
{\rm(1)} If the multi-saddle connection $D$ contains at least two multi-saddles at least one of which is outside of the boudnary, then there is a smooth Hamiltonian vector field $X' \in \mathcal{H}^r_*(i,g,p)$ obtained by its arbitrarily small $C^\infty$-perturbation in $\mathcal{H}^r_*(i,g,p)$ with $\mathrm{codim}(X) - \mathrm{codim}(X') = 1 = N(X') - N(X) = M(X') - M(X)$ such that $X'$ has same type of singular points of $X$ as in the top on Figure~\ref{splitting_02}. 
\\
{\rm(2)} If the multi-saddle connection $D$ contains a $\partial$-$k$-saddle {\rm(}$k \in \Z_{\geq 1}${\rm)}, then there is a smooth vector field $X'' \in \mathcal{H}^r_*(i,g,p)$ obtained by its arbitrarily small $C^\infty$-perturbation in $\mathcal{H}^r_*(i,g,p)$ with $\mathrm{codim}(X) - \mathrm{codim}(X'') = 1 = M(X'') - M(X)$ as in the second from bottom on Figure~\ref{splitting_02} such that the support of $X - X''$ is an isolated \nbd of the $\partial$-$k$-saddle of $X$. 
\end{lemma}

A Hamiltonian vector field $Y \in \mathcal{H}_*(i,g,p)$ is {\bf obtained by splitting of a multi-saddle connection} from a Hamiltonian vector field $X$ in $\mathcal{H}_*(i,g,p)$ if $Y$ is one of the following: 
\\
{\rm(1)} $Y$ has same type of singular points of $X$, the number of multi-saddle connections of $Y$ is exactly one more than one of $X$ (i.e. $N(Y) - N(X) = 1$), and for any \nbd $\mathcal{U} \subset \mathcal{H}_*(i,g,p)$ of $X$ there is $Y' \in \mathcal{U}$ which is topologically equivalent to $Y$. 
\\
{\rm(2)} $Y$ is obtained by pinching of multi-saddles on a same boundary of $X$, and for any \nbd $\mathcal{U} \subset \mathcal{H}_*(i,g,p)$ of $X$ there is $Y' \in \mathcal{U}$ which is topologically equivalent to $Y$. 

By definition, we obtain that $\mathop{\mathrm{codim}}(X) - \mathop{\mathrm{codim}}(Y) =1$ for any Hamiltonian vector field $Y$ obtained by splitting of a multi-saddle connection from a Hamiltonian vector field $X$ in $\mathcal{H}_*(i,g,p)$. 
Notice that $X', X'', X'''$ in Lemma~\ref{lem:perturbation_04} and the previous lemma are obtained by splitting of a multi-saddle connection from $X$. 
We have the following key lemma. 

\begin{lemma}\label{lem:open_dense02}
For 
any integers $k, i,g,p \in \Z_{\geq 0}$, the subset $\mathcal{H}_{*,k}(i,g,p)$ is open dense in the subspace $\mathcal{H}_{*, >k-1}(i,g,p)$ and any vector field in the subspace $\mathcal{H}_{*,k}(i,g,p)$ is $C^r$-structurally stable in the subspace $\mathcal{H}_{*, >k-1}(i,g,p)$. 
\end{lemma}

\begin{proof}
The openness and $C^r$-structural stability follow from Corollary~\ref{lem:open}. 

We claim that $\mathcal{H}_{*,k}(i,g,p)$ is dense. 
Indeed, fix a vector field $X'$ in $\mathcal{H}_{*, >k-1}(i)$.
By Lemma~\ref{lem:perturbation_03}, we can apply the Whitehead moves in $\mathcal{H}_{*, >k-1}(i,g,p)$ as possible. 
By Lemma~\ref{lem:perturbation_04}, using bump functions each of whose supports is near singular points or a boundary component, we can apply breaking operations of non-self-connected separatrices in $\mathcal{H}_{*, >k-1}(i,g,p)$ by changing the values of singular points for the Hamiltonians. 
Then the resulting vector field of $X'$ by a small perturbation in $\mathcal{H}_{*, >k-1}(i,g,p)$ belongs to $\mathcal{H}_{*,k}(i,g,p)$. 
\end{proof}

\subsubsection{Hamiltonian vector fields obtained by splitting of a multi-saddle connection}

We characterize the codimension one property as follows.

\begin{theorem}\label{lem:perturbation_05}
For any $X \in \mathcal{H}_*(i,g,p)$, there is its \nbd $\mathcal{U}$ in $\mathcal{H}_*(i,g,p)$ such that the following statements are equivalent for any $Y \in \mathcal{U}$:
\\
{\rm(1)} $\mathrm{codim}(X) - \mathrm{codim}(Y) =1$. 
\\
{\rm(2)} The Hamiltonian vector field $Y$ is obtained from $X$ by either a Whitehead move or  splitting of a multi-saddle connection. 
\end{theorem}

\begin{proof}
Take a \nbd $\mathcal{U}$ in $\mathcal{H}_*(i,g,p)$ of $X$ as in Lemma~\ref{lem:no_merge_msc},  Lemma~\ref{lem:perturbation}, and Proposition~\ref{prop:perturbation}. 
Fix any $Y \in \mathcal{U}$. 
By Lemma~\ref{lem:perturbation_03} and definition of splitting of a multi-saddle connection, 
the assertion {\rm(2)} implies the assertion {\rm(1)}. 

Suppose that $\mathrm{codim}(X) - \mathrm{codim}(Y) =1$. 
We claim that if $Y$ has same type of multi-saddles of $X$ then $Y$ is obtained by splitting of a multi-saddle connection. 
Indeed, suppose that $Y$ has same type of multi-saddles of $X$. 
Then $Y$ has same type of singular points of $X$. 
This implies that $n_{\partial}(Y) = n_{\partial}(Y)$ and  $n_{m}(Y) = n_{m}(X)$. 
Therefore $M(Y) - M(X) = 1$ because $\mathrm{codim}(X) - \mathrm{codim}(Y) =1$. 
By construction, we have that $N_{\mathrm{per}}(X) = N_{\mathrm{per}}(Y)$ and so that $N(X) + 1 = N(Y)$. 
By Lemma~\ref{lem:no_merge_msc}, since $Y$ has same type of singular points of $X$, 
the multi-saddle connection diagram $D(Y)$ is obtained by gluing the disks which are \nbds of isolated multi-saddles as in Figure~\ref{multi-saddle_nbd_std} along transverse closed arcs $I_{t,j}$ as in Lemma~\ref{lem:no_merge_msc} such that the gluing pairs of transverse closed arcs of $Y$ are the same as one of $X$ but the gluing points in the arcs are different. 
This means that $Y$ is obtained from $X$ by splitting of a multi-saddle connection. 

Thus we may assume that $Y$ does not have same type of multi-saddles of $X$. 
Then there is an isolated \nbd $V_i$ of a multi-saddle $x_i$ of $X$ as in Lemma~\ref{lem:no_merge_msc} such that $Y|_{V_i}$ contains multi-saddles $x'_{i,1}, \ldots , x'_{i,s_i}$ ($s_i \geq 2$). 

We claim that we may assume that $N_{\partial}(Y) = N_{\partial}(X)$. 
Indeed, assume that $N_{\partial}(Y) > N_{\partial}(X)$. 
By $N(Y) \geq N(X)$ and $n(Y) \geq n(X)$, from definition of $M(Z) = N_{\mathrm{per}}(Z) + N(Z)$ for any vector field $Z \in \mathcal{H}_*(i,g,p)$, Lemma~\ref{lem:perturbation} and Corollary~\ref{cor:deg_codim} and the codimension one property imply that $n(Y) = n(X)$, $N(Y) = N(X)$, and $N_{\mathrm{per}}(Y) = N_{\mathrm{per}}(X) + 1$. 
Then there is a multi-saddle connection $D$ of $X$ which intersects the boundary and whose continuation $D_Y$ for $Y$ contains a connected component which does not intersect the boundary. 
Lemma~\ref{lem:deg_codim} implies that $\mathrm{codim}(D) - \mathrm{codim}(D_Y) = 1$ and that every multi-saddle connection except $D$ of $X$ is isomorphic to the continuation as a surface graph. 
From $n(Y) = n(X)$, by defintion of pinching of multi-saddles on a same boundary, Lemma~\ref{lem:inv_index} implies that $Y$ is obtained from $X$ by pinching of multi-saddles on a same boundary. 

We claim that we may assume that $N(Y) = N(X)$. 
Indeed, assume that $N(Y) > N(X)$. 
By $N_{\mathrm{per}}(Y) \geq N_{\mathrm{per}}(X)$ and $n(Y) \geq n(X)$, from definition of $M(Z) = N_{\mathrm{per}}(Z) + N(Z)$ for any vector field $Z \in \mathcal{H}_*(i,g,p)$, Lemma~\ref{lem:perturbation} and Corollary~\ref{cor:deg_codim} and the codimension one property imply that $n(Y) = n(X)$, $N(Y) = N(X) + 1$, and $N_{\mathrm{per}}(Y) = N_{\mathrm{per}}(X)$. 
Then there is a multi-saddle connection $D$ of $X$ whose continuation $D_Y$ for $Y$ is not connected. 
Lemma~\ref{lem:deg_codim} implies that $\mathrm{codim}(D) - \mathrm{codim}(D_Y) = 1$ and that every multi-saddle connection except $D$ of $X$ is isomorphic to the continuation as a surface graph. 
From $n(Y) = n(X)$, by defintion of splitting of a multi-saddle connection, Lemma~\ref{lem:inv_index} implies that $Y$ is obtained from $X$ by splitting of a multi-saddle connection. 

Thus $M(Y) = M(X)$. 
The codimension one property implies that $n(Y) - n(X) = 1$ and so that $n_{\partial}(Y) + n_{m}(Y) = n_{\partial}(X) + n_{m}(X) + 1$. 
Then $Y|_{V_i}$ contains exactly two multi-saddles $x'_{i,1}$ and $x'_{i,2}$ for some $i$ and other restrictions $Y|_{V_j}$ ($j\neq i$) contains exactly one multi-saddles. 
By $M(Y) = M(X)$, the multi-saddles $x'_{i,1}$ and $x'_{i,2}$ are contained in a same multi-saddle connection of $Y$. 
From the Poincar{\'e}-Hopf theorem and the invariance of the number of centers, the restriction to $V_i$ of the multi-saddle connection of $Y$ containing $x'_{i,1}$ and $x'_{i,2}$ is a tree as in the right on Figure~\ref{splitting}. 
Therefore the multi-saddles $x'_{i,1}$ and $x'_{i,2}$ are obtained by a Whitehead move. 
The proof of the last claim in the proof of Lemma~\ref{lem:deg_codim} implies that muti-saddle connections $D_{Y,i}$ is obtained by a Whitehead move from $D_{X,i}$. 
The last claim in the proof of Lemma~\ref{lem:deg_codim} implies that muti-saddle connections $D_j(Y)$ and $D_j(X)$ for any $j \neq i$ are isomorphic as a surface graph. 
Therefore $Y$ is obtained from $X$ by a Whitehead move. 
\end{proof}


\section{Abstract cell complex structure and filtration of $\mathcal{H}_*(i)$}

The compactness of the surface implies the following finite property. 

\begin{lemma}\label{lem:top_eq_class}
For any $g, i, p \in \mathbb{Z}_{\geq 0}$ and any $r \in \Z_{\geq 0} \sqcup \{ \infty \}$, the union $\mathcal{H}^r_*(i,g,p)$ contains at most finitely many topological equivalence classes.  
\end{lemma}

\begin{proof}
Because the number of centers is bounded, so is the number of multi-saddles. 
Therefore there are at most finitely many multi-saddle connections that appear in $\mathcal{H}^r_*(i,g,p)$. 
This implies the finite possible combination. 
\end{proof}

\subsection{Filtration of open dense subspaces}

Corollary~\ref{lem:open}, Lemma~\ref{lem:open_dense02}, and Lemma~\ref{lem:top_eq_class} imply the following filtration of the space of Hamiltonian vector fields on compact surfaces. 

\begin{theorem}\label{main:01}
The following statements hold for  any $i, k, g, p \in \Z_{\geq 0}$ and any $r \in \Z_{\geq 0} \sqcup \{ \infty \}$:
\\
{\rm(1)} $\mathcal{H}^r_*(i,g,p) = \bigsqcup_{k=0}^\infty \mathcal{H}^r_{*,k}(i,g,p)$. 
\\
{\rm(2)} $\mathcal{H}^r_*(i,g,p) = \mathcal{H}^r_{*, >-1}(i,g,p) \supset \mathcal{H}^r_{*, >0}(i,g,p) \supset \mathcal{H}^r_{*, >1}(i,g,p) \supset \cdots$. 
\\
{\rm(3)} The subset $\mathcal{H}^r_{*,k}(i,g,p)$ is open in $\mathcal{H}^r_{*, >k-1}(i,g,p)$ and consists of $C^r$-structurally stable Hamiltonian vector fields in $\mathcal{H}^r_{*, >k-1}(i,g,p)$. 
\\
{\rm(4)} The subspace $\mathcal{H}^r_{*,k}(i,g,p)$ is $C^0$-dense in $\mathcal{H}^r_{*, >k-1}(i,g,p)$.
\end{theorem}

The author would like to know whether the subspace $\mathcal{H}^r_{*,k}(i,g,p)$ is $C^s$-dense in $\mathcal{H}^r_{*, >k-1}(i,g,p)$ for any $s \leq r \in \Z_{\geq 0} \sqcup \{ \infty \}$.

\subsubsection{Abstract cell complex structure}

For any $i, k, g, p \in \Z_{\geq 0}$ and any $r \in \Z_{\geq 0} \sqcup \{ \infty \}$, denote by $[\mathcal{H}^r_{*}(i,g,p)]$ the quotient space of $\mathcal{H}^r_{*}(i, g, p)$ by the topologically equivalence classes in $\mathcal{H}^r(g, p)$, and by $[\mathcal{H}^r_{*}(i, g, p)]$ the quotient space of $\mathcal{H}^r_{*}(i, g, p)$ by the topologically equivalence classes in $\mathcal{H}^r(g, p)$. 
Lemma~\ref{lem:perturbation_04+} implies the following statement. 

\begin{lemma}\label{lem:perturbation_04++}
For any $i, k, g, p \in \Z_{\geq 0}$, we have $\overline{[\mathcal{H}^\infty_{*,k}(i,g,p)]} = [\mathcal{H}^\infty_{*, >k-1}(i,g,p)]$. 
\end{lemma}

Lemma~\ref{lem:top_eq_class} implies that the quotient space $[\mathcal{H}^r_{*}(i, g, p)]$ is a finite $T_0$-space, and that the coheight of $[\mathcal{H}^r_{*}(i, g, p)]$ corresponds to the codimension as follows because of Theorem~\ref{main:01} and Lemma~\ref{lem:perturbation_04++}. 

\begin{theorem}\label{main:02}
The following statements hold for any $i, k, g, p \in \Z_{\geq 0}$ and any $r \in \Z_{\geq 0} \sqcup \{ \infty \}$:
\\
{\rm(1)} The subset $[\mathcal{H}^r_{*}(i, g, p)]$ is an abstract cell complex with respect to the specialization preorder and the codimension. 
\\
{\rm(2)} The quotient space $[\mathcal{H}^r_{*}(i, g, p)]$ is a finite $T_0$-space. 
\\
{\rm(3)} $[\mathcal{H}^r_{*}(i, g, p)]_{\dim [\mathcal{H}^r_{*}(i, g, p)]- k} = [\mathcal{H}^r_{*, k}(i, g, p)]$. 
\\
{\rm(4)} $\overline{[\mathcal{H}^r_{*}(i, g, p)]_k} = [\mathcal{H}^r_{*}(i, g, p)]_{\leq k} = [\mathcal{H}^r_{*, \geq k}(i, g, p)]$. 
\end{theorem}

This implies Theorem~\ref{main:02-}.

\subsection{Order complex of a finite $T_0$-space}

One can associate to any finite $T_0$-space $X$ the following simplicial complex, which is introduced by McCord \cite{mccord1966singular}, denoted by $\mathcal{K}(X)$, and called the {\bf order complex}: 
The set of vertices of $\mathcal{K}(X)$ is $X$ and the set of simplices is the set of chains on $X$. 
In \cite[Theorem~2]{mccord1966singular}, McCord proved that there is a weak homotopy equivalence from the geometric realization $\mathcal{K}(X)$ to $X$ (i.e. a continuous map $\mathcal{K}(X) \to X$ whose induced maps between all homotopy groups are isomorphisms). 
This implies the following simplicial complex structures of connected components of the space of topologically equivalence classes Hamiltonian vector field with finitely many singular points under the non-existence of creations and annihilations of singular points and fake multi-saddles. 


\begin{corollary}\label{cor:02}
The quotient space $[\mathcal{H}^r_{*}(i, g, p)]$ for any $i, g, p, r \in \Z_{\geq 0}$ is a simplicial complex up to weak homotopy equivalence. 
\end{corollary}

\section{Non-contractibility of the space of ``transitions of Hamiltonian vector fields on compact surfaces'}

\subsection{Topological properties of the spaces of Hamiltonian vector fields on compact surfaces}

We obtain the negative answer of the question of whether the structure of transitions of Hamiltonian vector fields on compact surfaces is contractible. 
More precisely, we have the following statement. 

\begin{theorem}\label{th:noncontractible}
The quotient space $[\mathcal{H}^r_{*}(3, 0, 1)]$ for any $r \in \Z_{\geq 0} \sqcup \{ \infty \}$ has non-contractible connected components. 
\end{theorem}

In fact, the quotient space $[\mathcal{H}^r_{*}(3, 0, 1)]$ has a connected component $X$ which is the weak homotopy type of a three-dimensional sphere. 
Denote by $\bm{[\mathcal{H}^r_{*}((i_-,i_+), g, p)]}$ the subspace of $[\mathcal{H}^r_{*}(i, g, p)]$ whose elements is the equivalence class of a vector fields such that $i_-$ {\rm(resp.} $i_+$) centers are surrounded by clockwise {\rm(resp.} counter-clockwise) periodic orbits. 
We will show that $X = [\mathcal{H}^r_{*}((1,2), 0, 1)]$ is the desired connected component. 

On the other hand, we have the following observation.

\begin{proposition}\label{prop:cont}
For any $i \in \Z_{>1}$ and any $r \in \Z_{\geq 0} \sqcup \{ \infty \}$, the connected component $[\mathcal{H}^r_{*}((1, i-1), 0, 0)]$ is contractible. 
\end{proposition}

To show the statements, we recall the theory of homotopy types of finite $T_0$-spaces as follows. 

\subsection{Homotopy types of finite $T_0$-spaces}

From now on, we equip a finite $T_0$-space $(X, \tau)$ with the specialization preorder $\leq_{\tau}$. 
Notice that $ x \leq_{\tau} y $ if and only if  $ x \in \overline{\{ y \}}$. 
Write the upset $\mathop{\uparrow}_{\tau} x  := \{ y \in X \mid x \leq_{\tau} y \}$ and the downset $\mathop{\downarrow}_{\tau} x  := \{ y \in X \mid y \leq_{\tau} x \}$. 
To state the characterization of contractibility, we recall beat points as follows \cite{may2003finite,stong1966finite}. 

\begin{definition}
A point $x \in X$ is a {\bf down beat point} {\rm(}or a {\bf colinear} point {\rm)} if there is a point $y \in X$ with  $\mathop{\downarrow}_{\tau} x - \{ x \} = \mathop{\downarrow}_{\tau} y$. 
\end{definition}

\begin{definition}
A point $x \in X$ is a {\bf up beat point} {\rm(}or a {\bf linear} point {\rm)} if there is a point $y \in X$ with  $\mathop{\uparrow}_{\tau} x - \{ x \} = \mathop{\uparrow}_{\tau} y$. 
\end{definition}

\begin{definition}
A point is a {\bf beat point} if it is a down or up beat point. 
\end{definition}

Notice that a point $x \in X$ is a beat point if and only if there is a point $y \in X$ such that either $\mathop{\uparrow}_{\tau} x - \{ x \} = \mathop{\uparrow}_{\tau} y$ or $\mathop{\downarrow}_{\tau} x - \{ x \} = \mathop{\downarrow}_{\tau} y$, and that the inclusion $X - \{x\} \to X$ for any beat point $x \in X$ is a strong deformation retract. 
A finite $T_0$-space without beat points is called a {\bf minimal finite space}.
A minimal finite space is a {\bf core} of a finite $T_0$-space $X$ if it is a strong deformation retract of $X$. 
In \cite[Theorem~4]{stong1966finite}, Stong proved that any finite $T_0$-space has the unique core up to homeomorphism and that two finite $T_0$-spaces are homotopy equivalent to each other if and only if their cores are homeomorphic to each other. 
In particular, a finite $T_0$-space $X$ is contractible if and only if there is a sequence $X = X_0 \supsetneq X_1 \supsetneq \cdots \supsetneq X_n = \{ * \}$, where $X_{i+1}$ is obtained from $X_i$ by removing a beat point.

\subsection{Contractibility of $[\mathcal{H}^r_{*}((1, i-1), 0, 0)]$}

We demonstrate the contractibility of some connected components as follows. 

\begin{proof}[Proof of Proposition~\ref{prop:cont}]

The equivalence class $[v]$ of a flow $v$ with $i$ centers and one $(i-2)$-multi-saddle on the sphere is the unique minimal element with respect to the specialization order. 
This means that any height one elements are beat points and so one can remove all height one elements of the resulting spaces obtained by removing beat points inductively. 
Therefore there is a sequence $[\mathcal{H}^r_{*}((1, i-1), 0, 0)] = X_0 \supsetneq X_1 \supsetneq \cdots \supsetneq X_n = \{ [v] \}$, where $n$ is the dimension of the abstract cell complex $[\mathcal{H}^r_{*}((1, i-1), 0, 0)]$. 
This means that $[\mathcal{H}^r_{*}((1, i-1), 0, 0)]$ is contractible. 
\end{proof}

\subsection{Non-contractibility of $[\mathcal{H}^r_{*}((1,2), 0, 1)]$}

From now on, fix any $r \in \Z_{\geq 0} \sqcup \{ \infty \}$. 
We have the following statement. 

\begin{proposition}\label{prop:5.3}
The quotient space $[\mathcal{H}^r_{*}((1,2), 0, 1)]$ is weakly homotopic to a three-dimensional sphere. 
\end{proposition}

The previous proposition implies This implies Theorem~\ref{main:01-} and Theorem~\ref{th:noncontractible}. 

A multi-saddle connection is a {\bf saddle connection} if any multi-saddles in it are either saddles or $\partial$-saddles. 
To demonstrate the previous proposition, we have the following statements.  

\begin{figure}
\begin{center}
\includegraphics[scale=0.375]{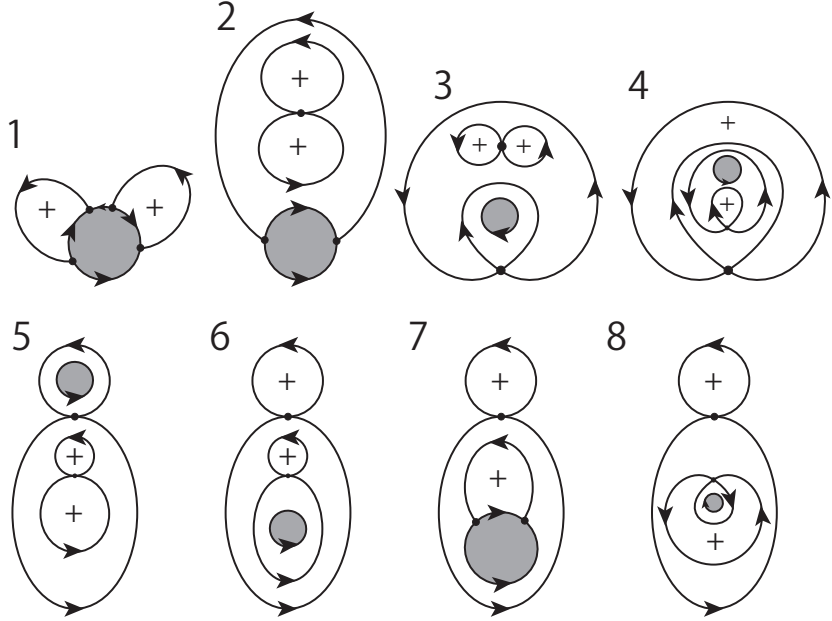}
\end{center}
\caption{Codimension zero elements in $[\mathcal{H}^r_{*}((1,2), 0, 1)]$.}
\label{fig:codimension_zero}
\end{figure}


\begin{lemma}\label{lem:001}
The subspace $[\mathcal{H}^r_{*,0}((1,2), 0, 1)]$ consists of eight topological equivalence class as in Figure~\ref{fig:codimension_zero}. 
\end{lemma}

\begin{proof}
Since there is only one center surrounded by clockwise periodic orbits, we can consider the center at the point at infinity and the topological equivalence classes of flows on the one-punctured sphere can be represented by multi-connection diagrams on the one-punctured plane surrounded by clockwise periodic orbits. 
Then we have exactly three possibilities of outermost saddle connections as on the upper of Figure~\ref{fig:codimension_zero}.
\begin{figure}
\begin{center}
\includegraphics[scale=0.225]{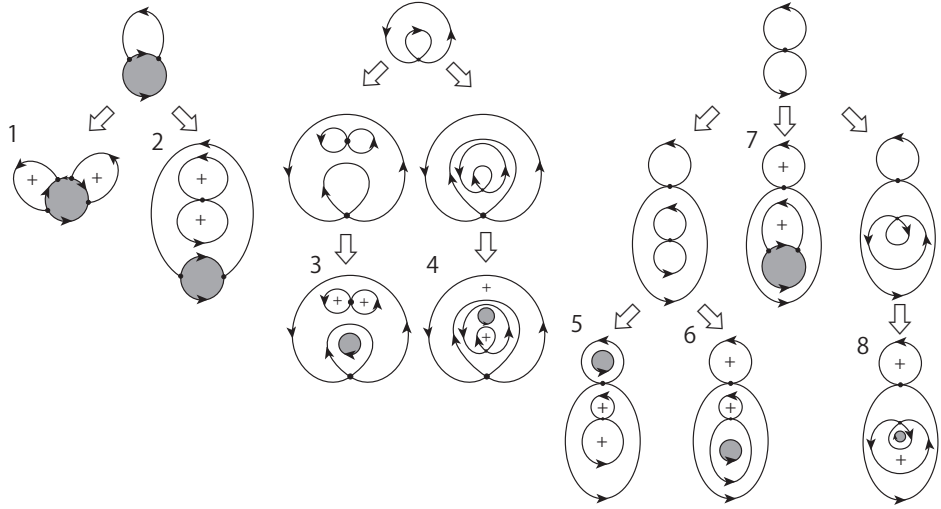}
\end{center}
\caption{The constructions of combinatorial structures of saddle connections of codimension zero elements.}
\label{fig:codimension_zero}
\end{figure} 
By definition of codimension, one can construct all possible codimension zero topological equivalence classes of flows by inserting saddle connections and periodic orbits on the boundary and by replacing saddle connections to add separatrices between $\partial$-saddles. 
Then the codimension zero topological equivalence classes of flows are represented by eight structurally stable Hamiltonian vector fields as in Figure~\ref{fig:codimension_zero}.  
\end{proof}

\begin{figure}
\begin{center}
\includegraphics[scale=0.225]{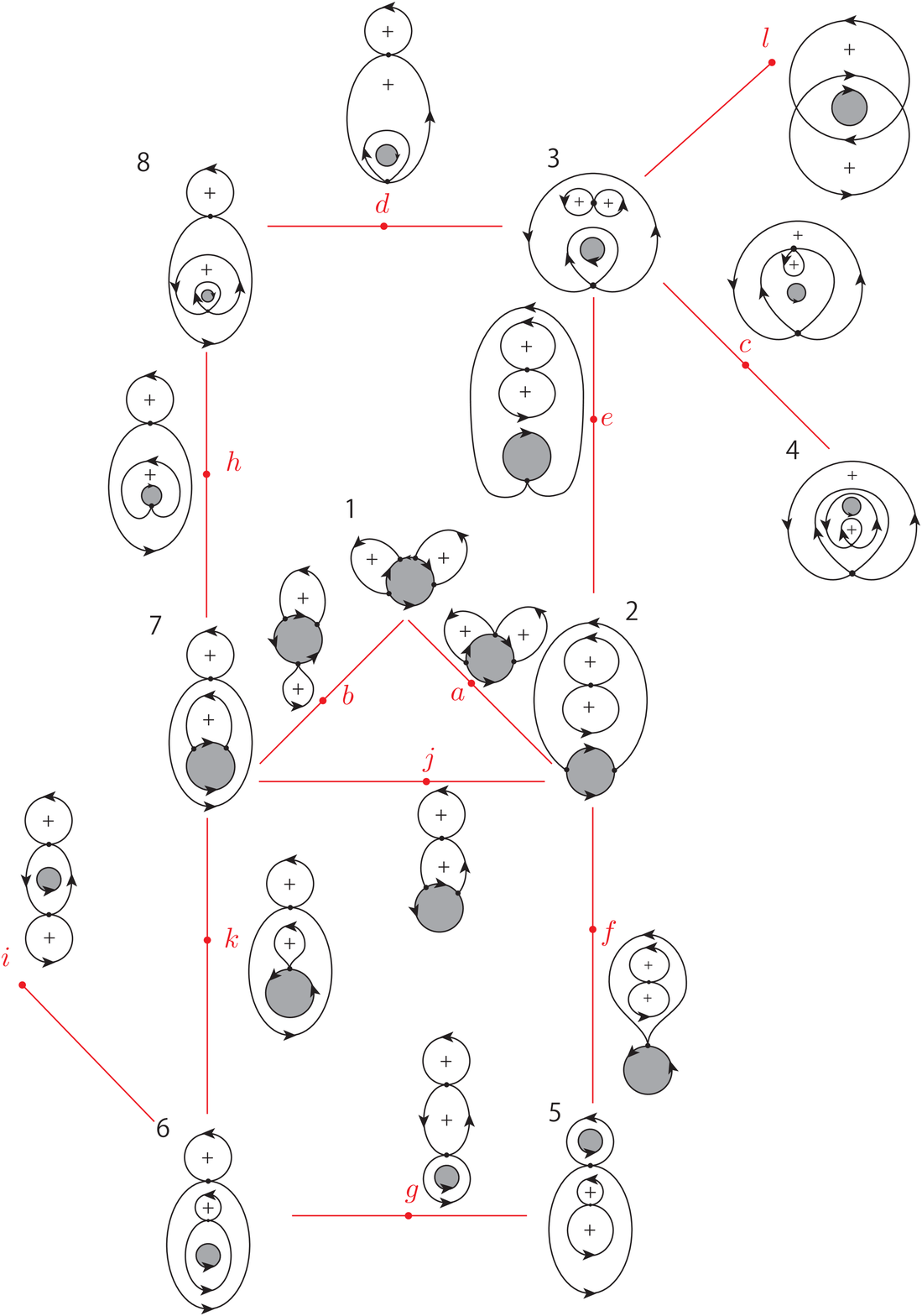}
\end{center}
\caption{Codimension one elements in $[\mathcal{H}^r_{*}((1,2), 0, 1)]$.}
\label{fig:codimension_one}
\end{figure}

\begin{lemma}\label{lem:001}
The subspace $[\mathcal{H}^r_{*,1}((1,2), 0, 1)]$ consists of twelve topological equivalence class as in Figure~\ref{fig:codimension_one}. 
\end{lemma}

\begin{proof}
Since any elements in $[\mathcal{H}^r_{*}((1,2), 0, 1)]$ has one boundary component, the codimension one elements have exactly one of the following structures: 
\\
{\rm(1)} A pair of two saddles connected by heteroclinic separatrices. 
\\
{\rm(2)} A triple of one saddle and two $\partial$-saddles connected by heteroclinic separatrices. 
\\
{\rm(3)} One $\partial$-$1$-saddle. 
\\
Therefore all codimension one elements with non-self-connected separatrices in $[\mathcal{H}^r_{*}((1,2), 0, 1)]$ are listed in Figure~\ref{fig:codimension_one_proof}. 
\begin{figure}
\begin{center}
\includegraphics[scale=0.15]{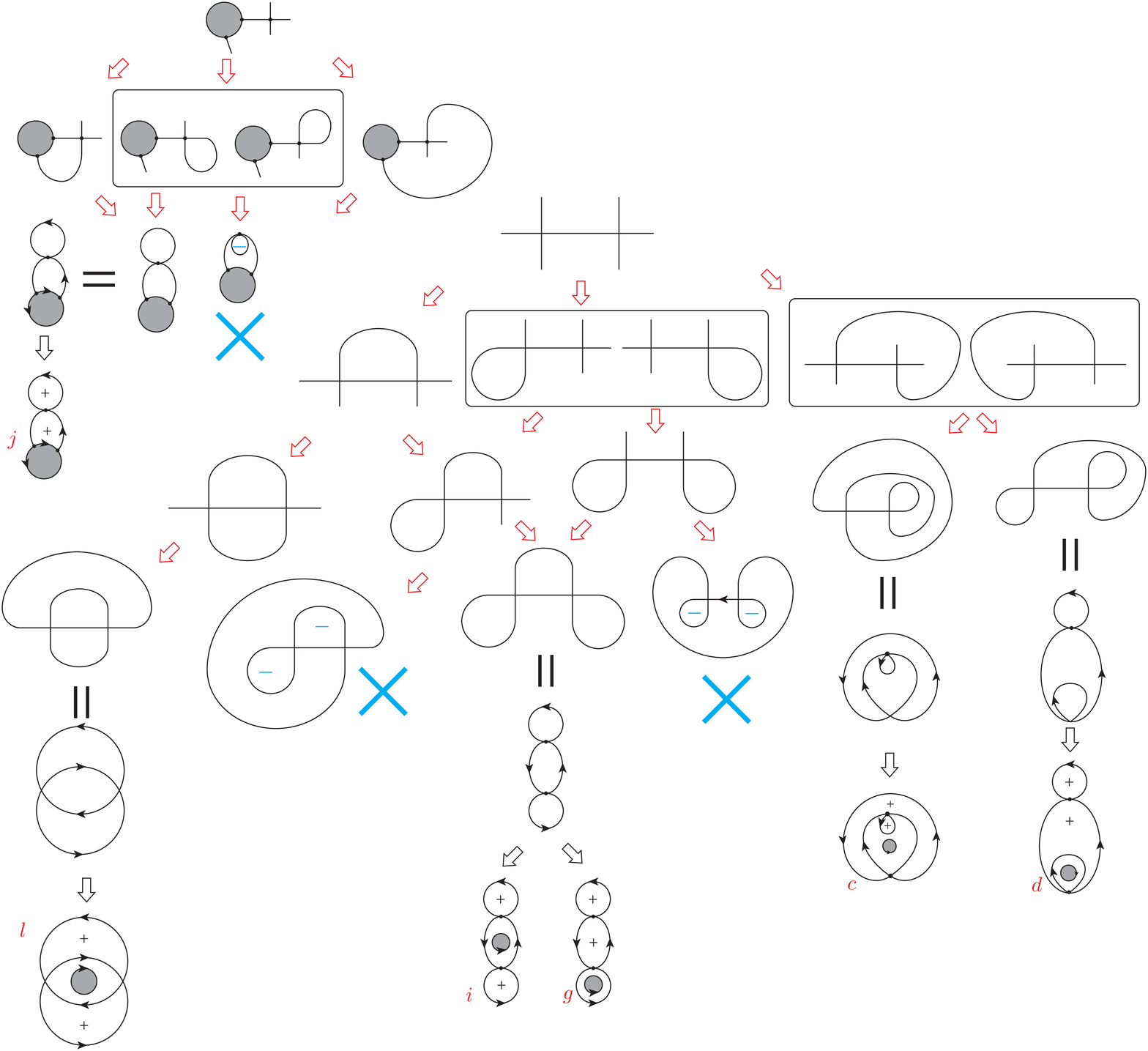}
\end{center}
\caption{The constructions of codimension one elements with non-self-connected separatrices in $[\mathcal{H}^r_{*}((1,2), 0, 1)]$.}
\label{fig:codimension_one_proof}
\end{figure} 
Moreover, all codimension one elements with one pinching in $[\mathcal{H}^r_{*}((1,2), 0, 1)]$ are listed in Figure~\ref{fig:codimension_one_proof_02}. 
\begin{figure}
\begin{center}
\includegraphics[scale=0.15]{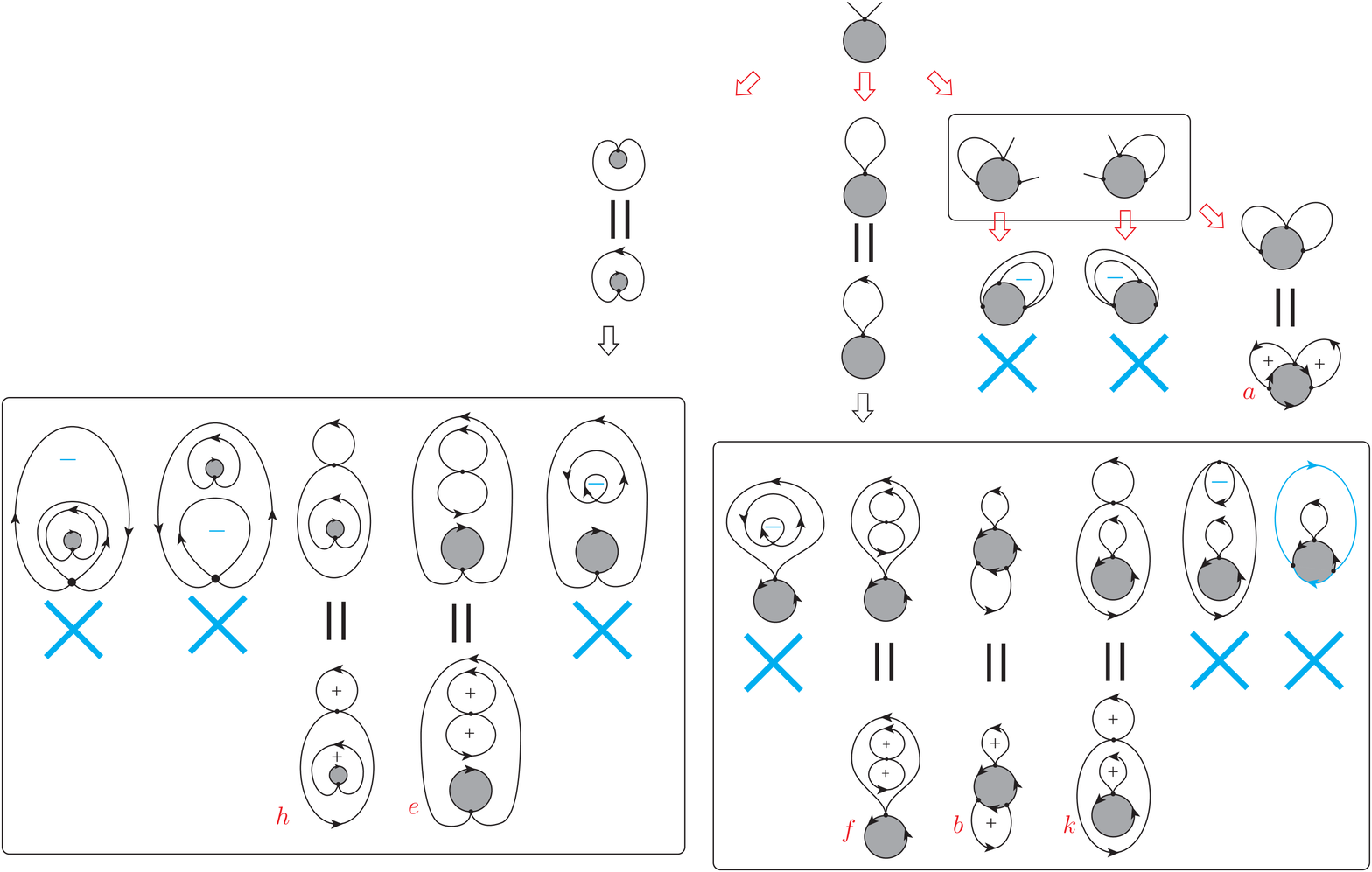}
\end{center}
\caption{The constructions of codimension one elements with one pinching in $[\mathcal{H}^r_{*}((1,2), 0, 1)]$.}
\label{fig:codimension_one_proof_02}
\end{figure} %
Then the list of twelve codimension one topological equivalence classes of flows are represented by six $h$-unstable Hamiltonian vector fields and six $p$-unstable Hamiltonian vector fields. 
\end{proof}

Note that we can also calculate the list of twelve codimension one topological equivalence classes of flows, checking the complete list of codimension one elements in \cite{yokoyama2021complete}. 
We list all codimension two elements in $[\mathcal{H}^r_{*}((1,2), 0, 1)]$.
\begin{figure}
\begin{center}
\includegraphics[scale=0.2]{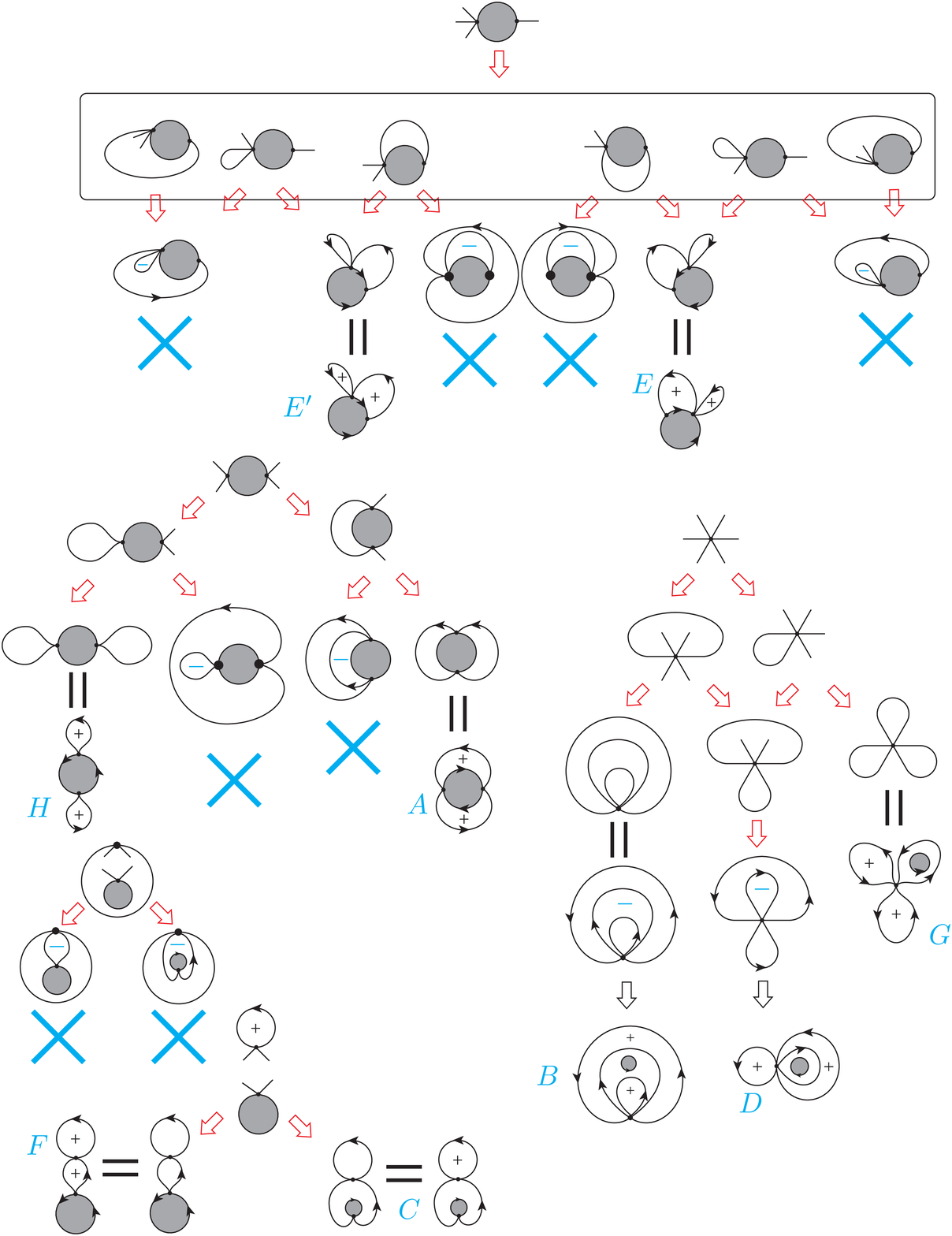}
\end{center}
\caption{The constructions of codimension two elements in $[\mathcal{H}^r_{*}((1,2), 0, 1)]$.}
\label{fig:codimension_two_proof}
\end{figure} %

\begin{figure}
\begin{center}
\includegraphics[scale=0.225]{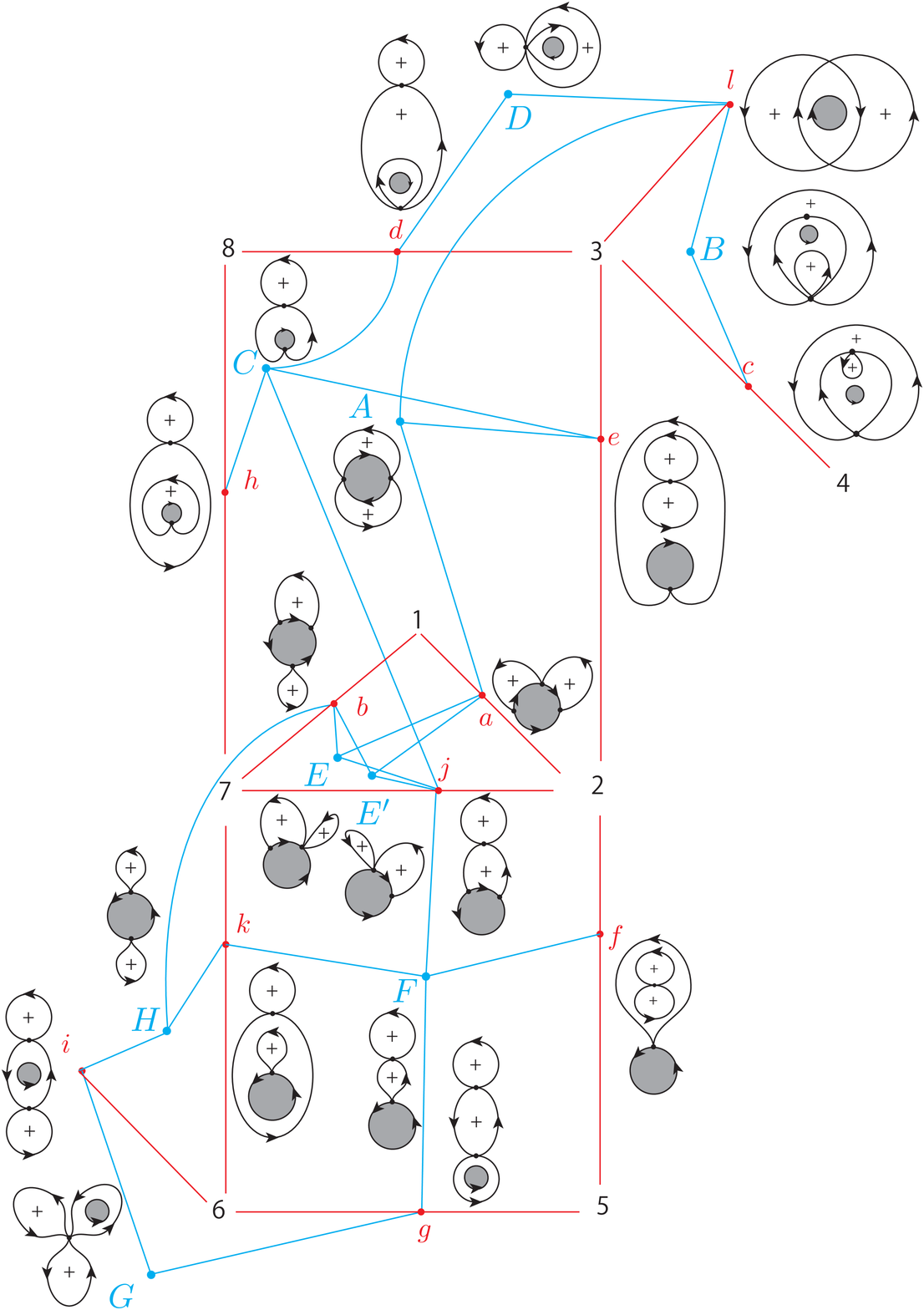}
\end{center}
\caption{Codimension two elements in $[\mathcal{H}^r_{*}((1,2), 0, 1)]$.}
\label{fig:codimension_two}
\end{figure} 

\begin{lemma}\label{lem:002}
The subspace $[\mathcal{H}^r_{*,2}((1,2), 0, 1)]$ consists of nine topological equivalence class as in Figure~\ref{fig:codimension_two}. 
\end{lemma}

\begin{proof}
By definition of codimension, any Hamiltonian vector fields whose topological equivalence classes are codimension two has one of the following structures as in Figure~\ref{fig:codimension_two_proof}: 
\\
{\rm(1)} Two $\partial$-$1$-saddles. 
\\
{\rm(2)} One $\partial$-$3/2$-saddle. 
\\
{\rm(3)} One $2$-saddle.
\\
{\rm(4)} One $\partial$-$1$-saddle with non-self-connected separatrices. 
\\
Merging two saddles into one $2$-saddle, we obtain seven multi-saddle connection diagrams $B$, $D$, and $G$ as on the middle right Figure~\ref{fig:codimension_two_proof}. 
Merging three $\partial$-saddles into one $\partial$-$3/2$-saddle, we obtain seven multi-saddle connection diagrams $E$ and $E'$ as on the upper of Figure~\ref{fig:codimension_two_proof}. 
Merging two pairs of two $\partial$-saddles into two $\partial$-$1$-saddles, we obtain seven multi-saddle connection diagrams $A$ and $H$ as on the middle left of  Figure~\ref{fig:codimension_two_proof}.  
Moreover, by merging a self-connected multi-saddle connection with a saddle and a multi-saddle connection with a pinching, we obtain two multi-saddle connection diagrams $C$ and $F$ as on the lower left of Figure~\ref{fig:codimension_two_proof}. 
\end{proof}

We list all codimension three elements in $[\mathcal{H}^r_{*}((1,2), 0, 1)]$.

\begin{figure}
\begin{center}
\includegraphics[scale=0.225]{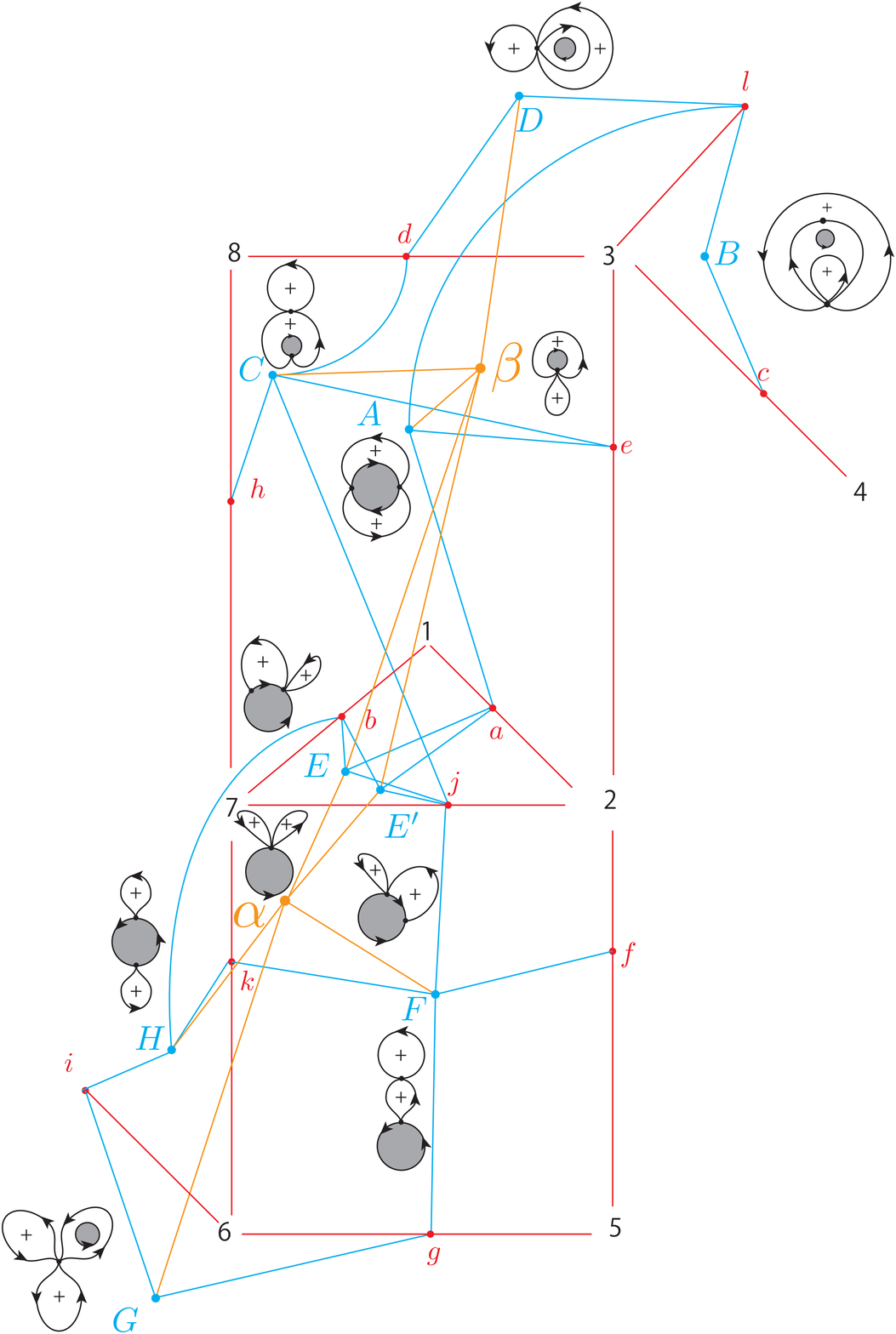}
\end{center}
\caption{Codimension three elements in $[\mathcal{H}^r_{*}((1,2), 0, 1)]$.}
\label{fig:codimension_three}
\end{figure}

\begin{lemma}\label{lem:003}
The subspace $[\mathcal{H}^r_{*,3}((1,2), 0, 1)]$ consists of two topological equivalence class as in Figure~\ref{fig:codimension_three}, each of which is represented by a flow with one $\partial$-$2$-saddle. 
\end{lemma}

\begin{proof}
By definition of codimension, any Hamiltonian vector fields whose topological equivalence classes are codimension three has exactly one $\partial$-$2$-saddle as in Figure~\ref{fig:codimension_three_proof}: 
By listing all the possible combinations, we obtain two multi-saddle connection diagrams $\alpha$ and $\beta$ as in Figure~\ref{fig:codimension_three_proof}. 
\begin{figure}
\begin{center}
\includegraphics[scale=0.2]{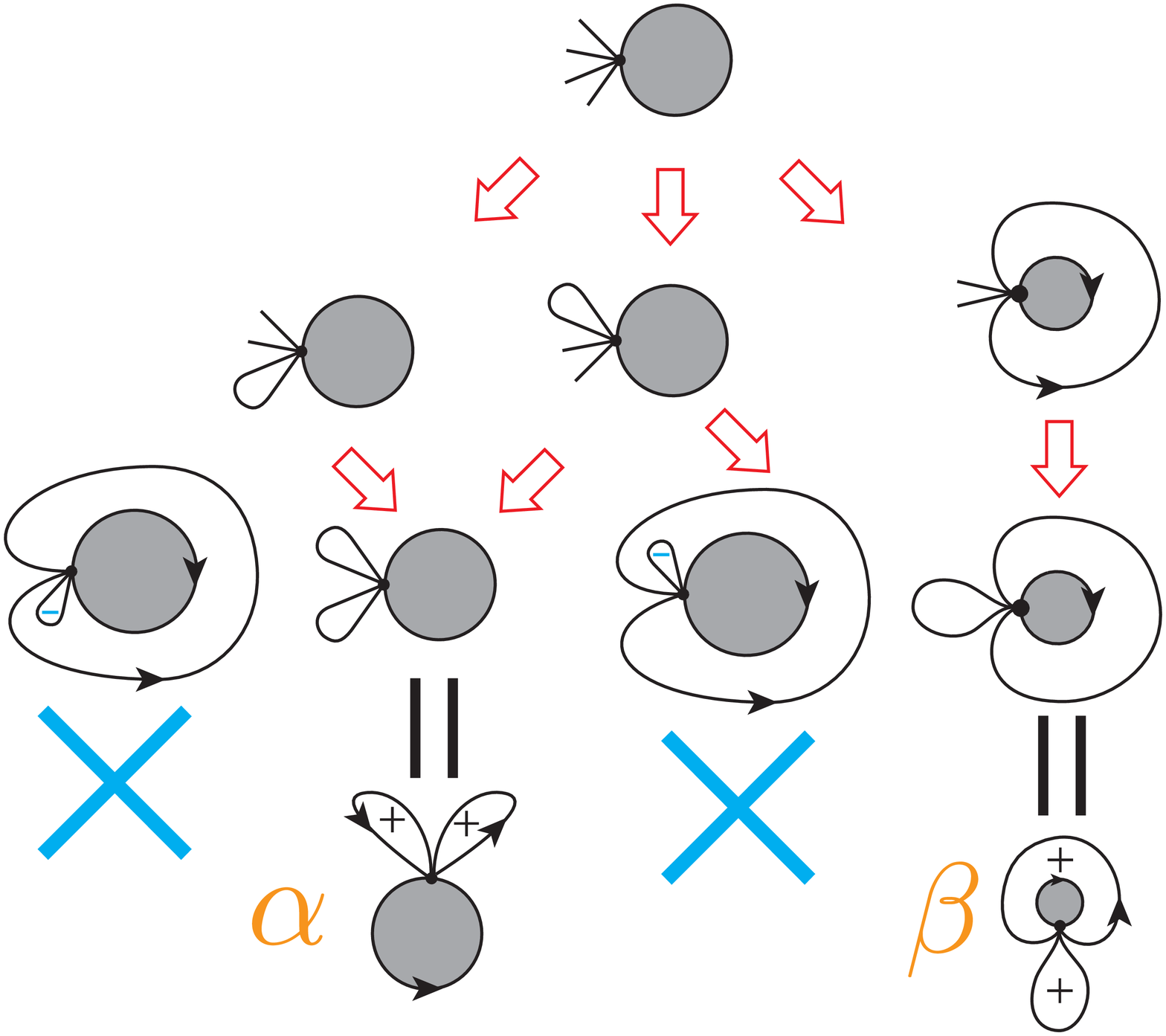}
\end{center}
\caption{The constructions of codimension three elements in $[\mathcal{H}^r_{*}((1,2), 0, 1)]$.}
\label{fig:codimension_three_proof}
\end{figure} 
\end{proof}

By the previous three lemmas, we have the following statement. 

\begin{figure}
\begin{center}
\includegraphics[scale=0.12]{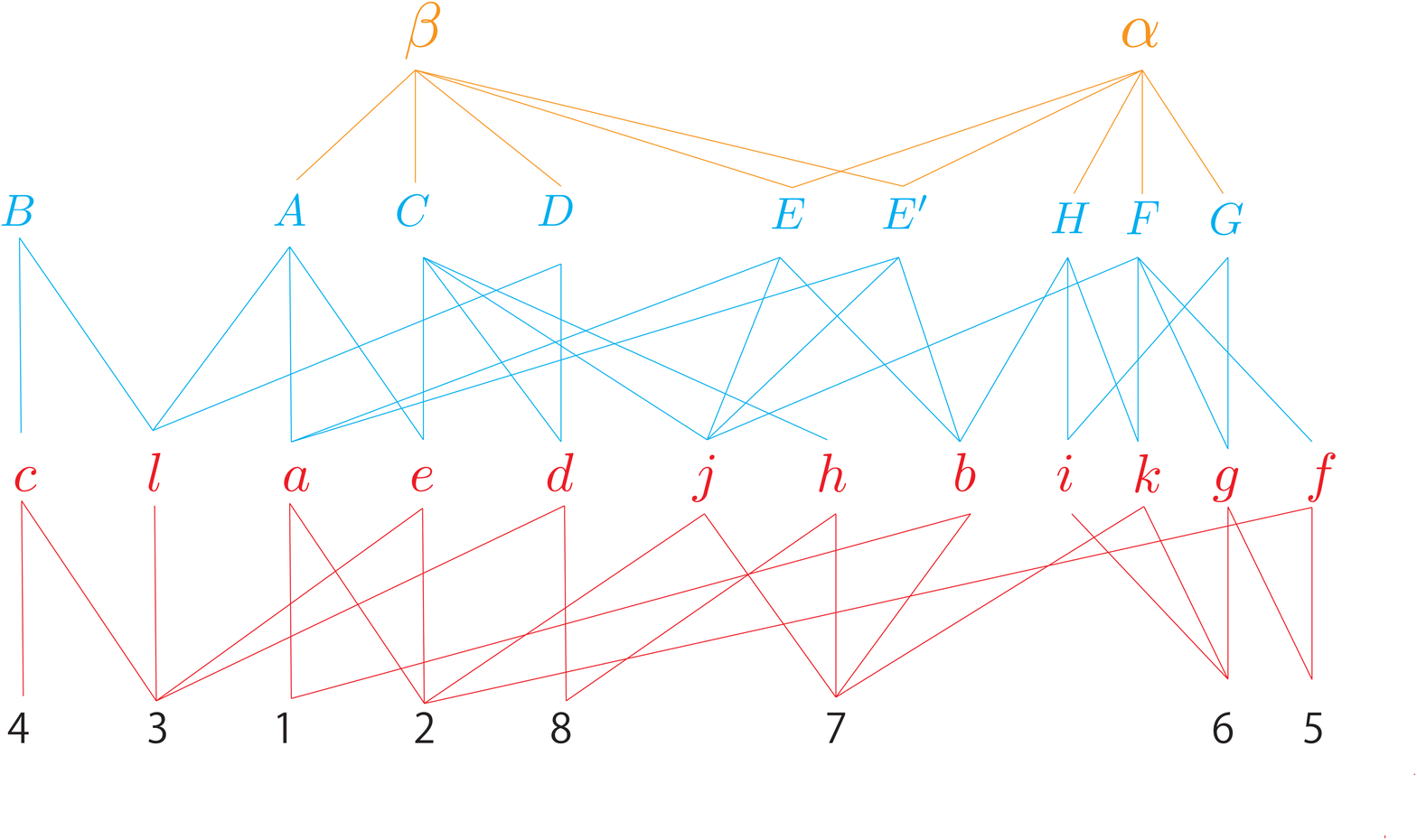}
\end{center}
\caption{Hesse diagram of the opposite order of the specialization order of the finite $T_0$-space $[\mathcal{H}^r_{*}((1,2), 0, 1)]$.}
\label{fig:Hessian}
\end{figure}

\begin{lemma}\label{lem:004}
The Hesse diagram of the opposite order of the specialization order of the finite $T_0$-space $[\mathcal{H}^r_{*}((1,2), 0, 1)]$ is shown in Figure~\ref{fig:Hessian}. 
\end{lemma}

\begin{figure}
\begin{center}
\includegraphics[scale=0.12]{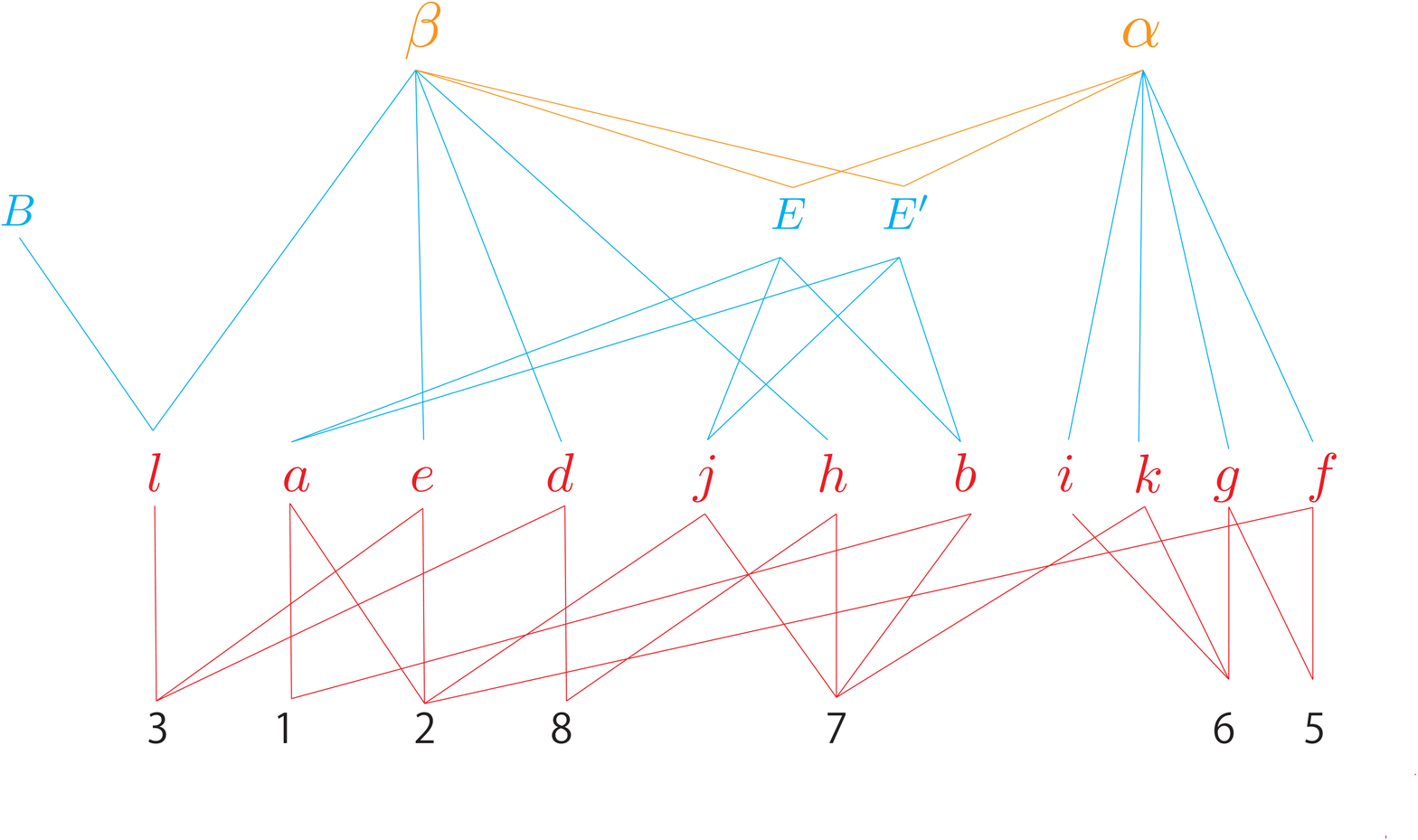}
\end{center}
\caption{The resulting Hesse diagram of Hesse diagram of the opposite order of the specialization order of the finite $T_0$-space $[\mathcal{H}^r_{*}((1,2), 0, 1)]$ by removing $A$, $C$, $D$, $F$, $G$, $H$, $4$, and $c$.}
\label{fig:Hessian_reduced_01}
\end{figure}

\begin{figure}
\begin{center}
\includegraphics[scale=0.12]{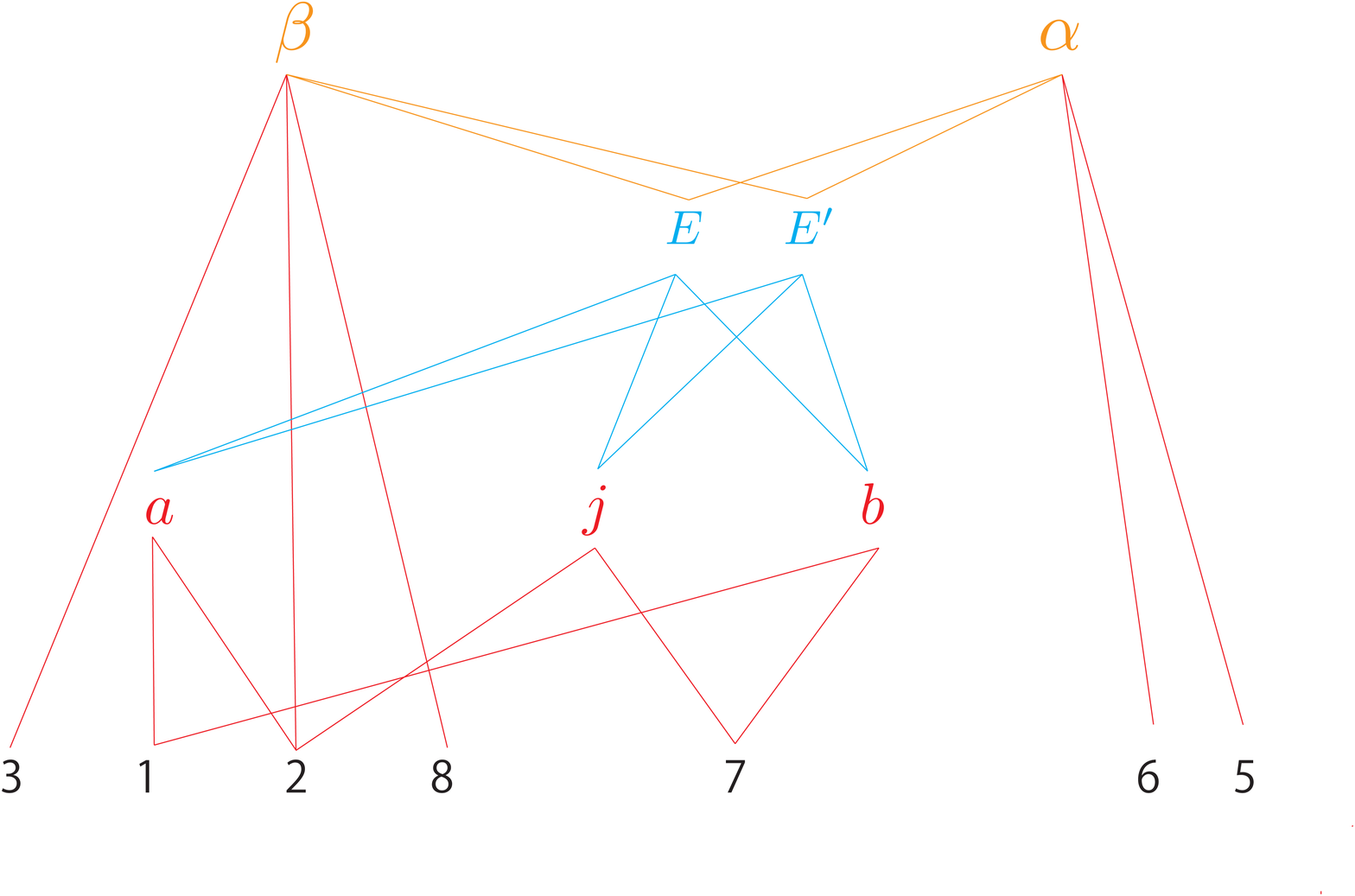}
\end{center}
\caption{The resulting Hesse diagram of Hesse diagram of the opposite order of the specialization order of the finite $T_0$-space $[\mathcal{H}^r_{*}((1,2), 0, 1)]$ by removing $A$, $C$, $D$, $F$, $G$, $H$, $4$, $c$, $B$, $l$, $e$, $d$, $h$, $i$, $k$, $g$, and $f$.}
\label{fig:Hessian_reduced_02}
\end{figure}

\begin{figure}
\begin{center}
\includegraphics[scale=0.12]{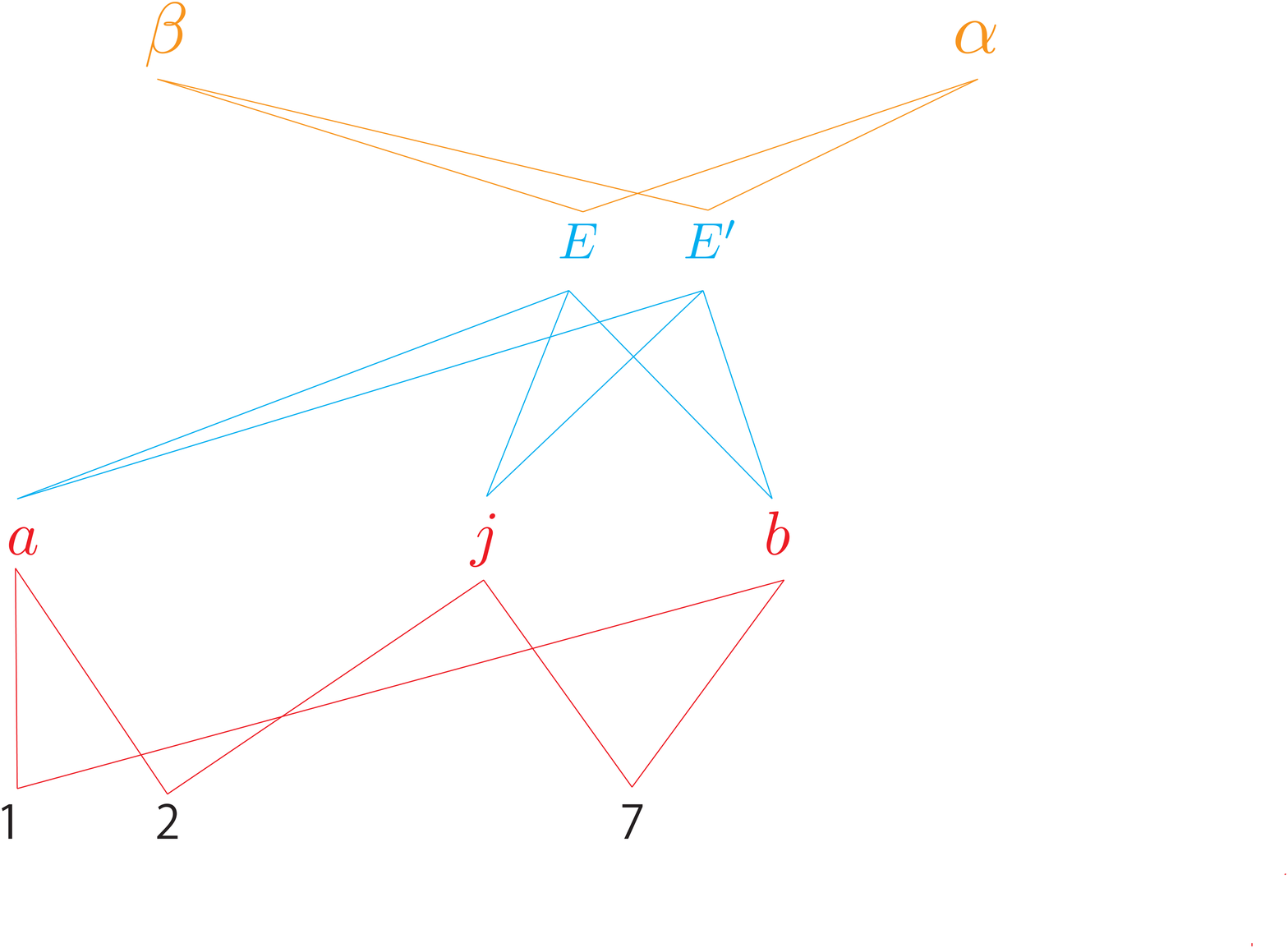}
\end{center}
\caption{Hesse diagram of the opposite order of the specialization order of the minimal finite space of the finite $T_0$-space $[\mathcal{H}^r_{*}((1,2), 0, 1)]$.}
\label{fig:Hessian_reduced_03}
\end{figure}

\begin{proof}
The Hesse diagram of the opposite order of the specialization order of the finite $T_0$-space $[\mathcal{H}^r_{*}((1,2), 0, 1)]$ is shown in Figure~\ref{fig:Hessian}.
Remove up beat points $A$, $C$, $D$, $F$, $G$, $H$, and $4$. 
Removing a down beat point $c$ recursively, we have the Hesse diagram as in Figure~\ref{fig:Hessian_reduced_01}. 
Remove down beat points $B$, $l$, and up down beat points $e$, $d$, $h$, $i$, $k$, $g$, and $f$ recursively, we obtain the Hesse diagram as in Figure~\ref{fig:Hessian_reduced_02}. 
Remove up beat points $3$, $8$, $6$, and $5$ recursively, we have the Hesse diagram of the opposite order of the specialization order of the minimal finite space of the finite $T_0$-space $[\mathcal{H}^r_{*}((1,2), 0, 1)]$ is shown in Figure~\ref{fig:Hessian_reduced_03}. 
Therefore $[\mathcal{H}^r_{*}((1,2), 0, 1)]$ is not contractible. 
Since the inclusion $X - \{x\} \to X$ into a finite $T_0$-space $X$ for any beat point $x \in X$ is a strong deformation retract, the minimal finite space of $[\mathcal{H}^r_{*}((1,2), 0, 1)]$ is homotopic to $[\mathcal{H}^r_{*}((1,2), 0, 1)]$. 
\end{proof}

\begin{lemma}\label{lem:006}
The minimal finite space of the finite $T_0$-space $[\mathcal{H}^r_{*}((1,2), 0, 1)]$ is weakly homotopic to the three-dimensional sphere. 
\end{lemma}

\begin{figure}
\begin{center}
\includegraphics[scale=0.2]{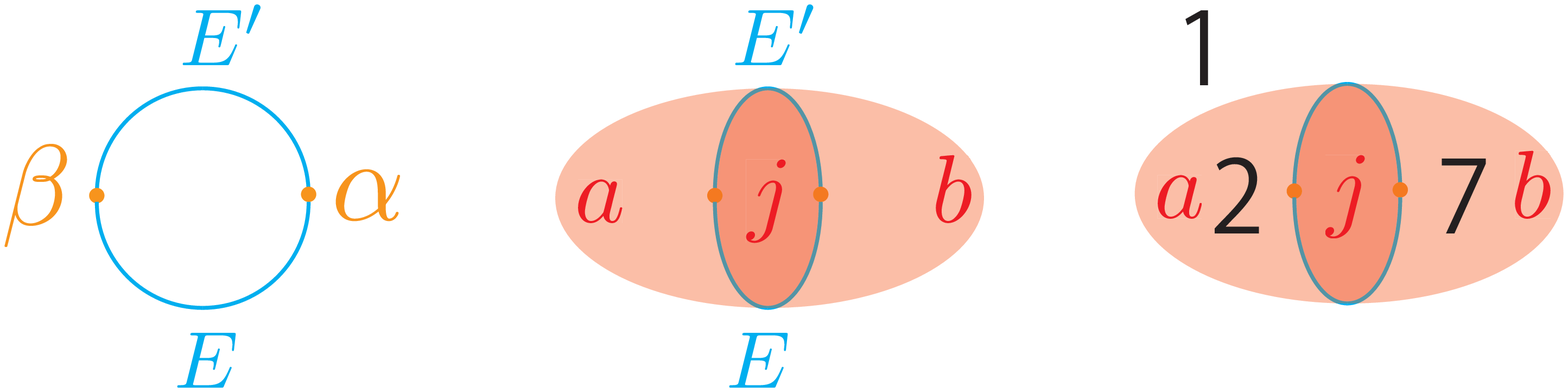}
\end{center}
\caption{Left, circle; middle, the resulting space of a closed three-dimensional ball by removing two open three-dimensional balls; right, a three-dimensional sphere.}
\label{fig:order_cpx01}
\end{figure}

\begin{proof}
The order complex $\mathcal{K}(X')$ of the subset $X' := \{ \alpha, \beta, E, E' \}$ is a circle as on the left of  Figure~\ref{fig:order_cpx01} and one of $\{ \alpha, \beta, E, E', a, b, j \}$ is the resulting space of a closed three-dimensional ball by removing two open three-dimensional balls as on the left of  Figure~\ref{fig:order_cpx01}, which is homotopic to a bouquet of two spheres. 
Since the order complex $\mathcal{K}(X'')$  of $X'' := \{ \alpha, \beta, E, E', a, b, j, 2, 7 \}$ is a closed three-dimensional ball, one of the subset $\{ \alpha, \beta, E, E', a, b, j, 1, 2, 7 \}$ is a three-dimensional sphere. 
By \cite[Theorem~2]{mccord1966singular}, the order complex $\mathcal{K}(X''')$  of the minimal finite space $X''' := \{ \alpha, \beta, E, E', a, b, j, 2, 7 \}$ of $[\mathcal{H}^r_{*}((1,2), 0, 1)]$ is weak homotopy equivalent to the minimal finite space $\{ \alpha, \beta, E, E', a, b, j, 2, 7 \}$. 
Because the inclusion $X - \{x\} \to X$ into a finite $T_0$-space for any beat point $x \in X$ is a strong deformation retract, the subspace $[\mathcal{H}^r_{*}((1,2), 0, 1)]$ is weak homotopy equivalent to a three-dimensional sphere. 
\end{proof}

The previous lemma implies Proposition~\ref{prop:5.3}. 
We could like to ask the following question. 
\begin{question}
Is every connected component $[\mathcal{H}^r_{*}((i_-,i_+), g, p)]$ simply connected? 
\end{question}

\section{Final remark}

In \cite{2022VladislavMorse}, it is shown that any connected components of the space of gradient vector fields on a compact surface have finite combinatorial structures under the non-existence of creations and annihilations of singular points. 
Our combinatorial method to calculate weak homotopy class also can be applied to the following question for the two-dimensional case. 

\begin{question}
Does the space of topological equivalence classes of gradient vector fields on a compact manifold have non-contractible connected components under the non-existence of creations and annihilations of singular points? 
\end{question}

In fact, it is shown that there are such connected components for the two-dimensional case \cite{yokoyama2022combinatorial}. 

\bibliographystyle{abbrv}
\bibliography{../yt20211011}

\end{document}